\documentclass[12pt]{article}
\usepackage{amssymb,latexsym,amsmath}
\usepackage{epsfig}
\usepackage{a4wide}
\usepackage{graphicx}

\def\dist{{\rm dist \, }}

\def\Re{{\rm Re \, }}

\def\ds{\displaystyle}

\newtheorem{theorem}{Theorem}[section]

\newtheorem{proposition}[theorem]{Proposition}

\newtheorem{Definition}[theorem]{Definition}
\newenvironment{definition}{\begin{Definition}\rm}{\end{Definition}}
\newtheorem{Remark}[theorem]{Remark}
\newenvironment{remark}{\begin{Remark}\rm}{\end{Remark}}
\newtheorem{Example}[theorem]{Example}

\newenvironment{proof}%
{\rm \trivlist \item[\hskip \labelsep{\bf Proof. }]}%
{\hspace*{\fill}$\Box$\endtrivlist}
\newenvironment{varproof}%
{\rm \trivlist \item[\hskip \labelsep{\bf Proof}]}%
{\hspace*{\fill}$\Box$\endtrivlist}

\begin{document}

\begin{center} \Large\bf
Asymptotic zero behavior of Laguerre polynomials with negative parameter
\end{center}

    \begin{center} \large
        A.B.J. Kuijlaars\footnote{Supported by FWO research project G.0176.02 and
        by INTAS project 00-272.} \\
        \normalsize \em
        Department of Mathematics, Katholieke Universiteit Leuven, \\
        Celestijnenlaan 200 B,  3001 Leuven, Belgium \\
        \rm arno@wis.kuleuven.ac.be \\[3ex]
        \rm and \\[3ex]
        \large
        K.T-R McLaughlin\footnote{Supported by NSF grant \#DMS-9970328.}\\
        \normalsize \em
        Department of Mathematics, University of North Carolina, \\
        Chapel Hill, NC 27599, U.S.A. \\
        \rm and \\
        \em
        Department of Mathematics, University of Arizona, \\
        Tucson, AZ 85721, U.S.A.\\
        \rm mcl@amath.unc.edu \\[10pt]
    \end{center}

\begin{center}
\end{center}

\begin{abstract}
We consider Laguerre polynomials $L_n^{(\alpha_n)}(nz)$ with varying
negative parameters $\alpha_n$, such that the limit $A = -\lim_n \alpha_n/n$
exists and belongs to $(0,1)$. For $A > 1$, it is known that the zeros
accumulate along an open contour in the complex plane.
For every $A \in (0,1)$, we describe a one-parameter family of possible
limit sets of the zeros.
Under the condition that the limit $r= - \lim_n \frac{1}{n} \log
\left[\dist(\alpha_n, \mathbb Z)\right]$ exists, we show that
the zeros accumulate on $\Gamma_r \cup [\beta_1,\beta_2]$ with
$\beta_1$ and $\beta_2$ only depending on $A$. For $r \in [0,\infty)$,
$\Gamma_r$ is a closed loop encircling the origin, which for $r = +\infty$,
reduces to the origin. This shows a great sensitivity of the
zeros to $\alpha_n$'s proximity to the integers.

We use a Riemann-Hilbert formulation for the Laguerre polynomials,
together with the steepest descent method of Deift and Zhou
to obtain asymptotics for the polynomials, from which the zero behavior follows.
\end{abstract}



\section{Introduction}
The zeros of Laguerre polynomials $L_n^{(\alpha)}$ with negative
parameters $\alpha$ have  a remarkable behavior. For $\alpha > -1$,
these polynomials are orthogonal with respect to the weight $x^{\alpha} e^{-x}$
on $[0,\infty)$, and the zeros are all real and positive.
This is no longer true if $\alpha < -1$.
The explicit representation, see \cite{Szego},
\begin{equation} \label{eq11}
    L_n^{(\alpha)}(z)  = \sum_{k=0}^n \binom{n+ \alpha}{n-k}
    \frac{(-z)^k}{k!}
\end{equation}
shows that $L_n^{(\alpha)}$ has a zero of order $|\alpha|$ at $0$ if $\alpha$
is an integer in $\{-n,-n+1, \ldots, -1\}$.  For these values of $\alpha$ we have
\begin{equation} \label{eq11a}
    L_n^{(\alpha)}(z) = \frac{(n+\alpha)!}{n!}
    (-z)^{-\alpha} L_{n+\alpha}^{(-\alpha)}(z).
\end{equation}
This means that in addition to the $|\alpha|$ zeros at $0$, there are $n- |\alpha|$
zeros in $(0,\infty)$.
For a non-integer $\alpha \in (-n,-1)$ it is known that there are
exactly $n - [-\alpha]$ positive zeros and for
$\alpha < -n$ there are no positive zeros at all, see \cite[\S 6.73]{Szego}.
In all cases there is at most one negative real zero.

\begin{figure}[th!]
\centerline{\includegraphics[width=5cm]{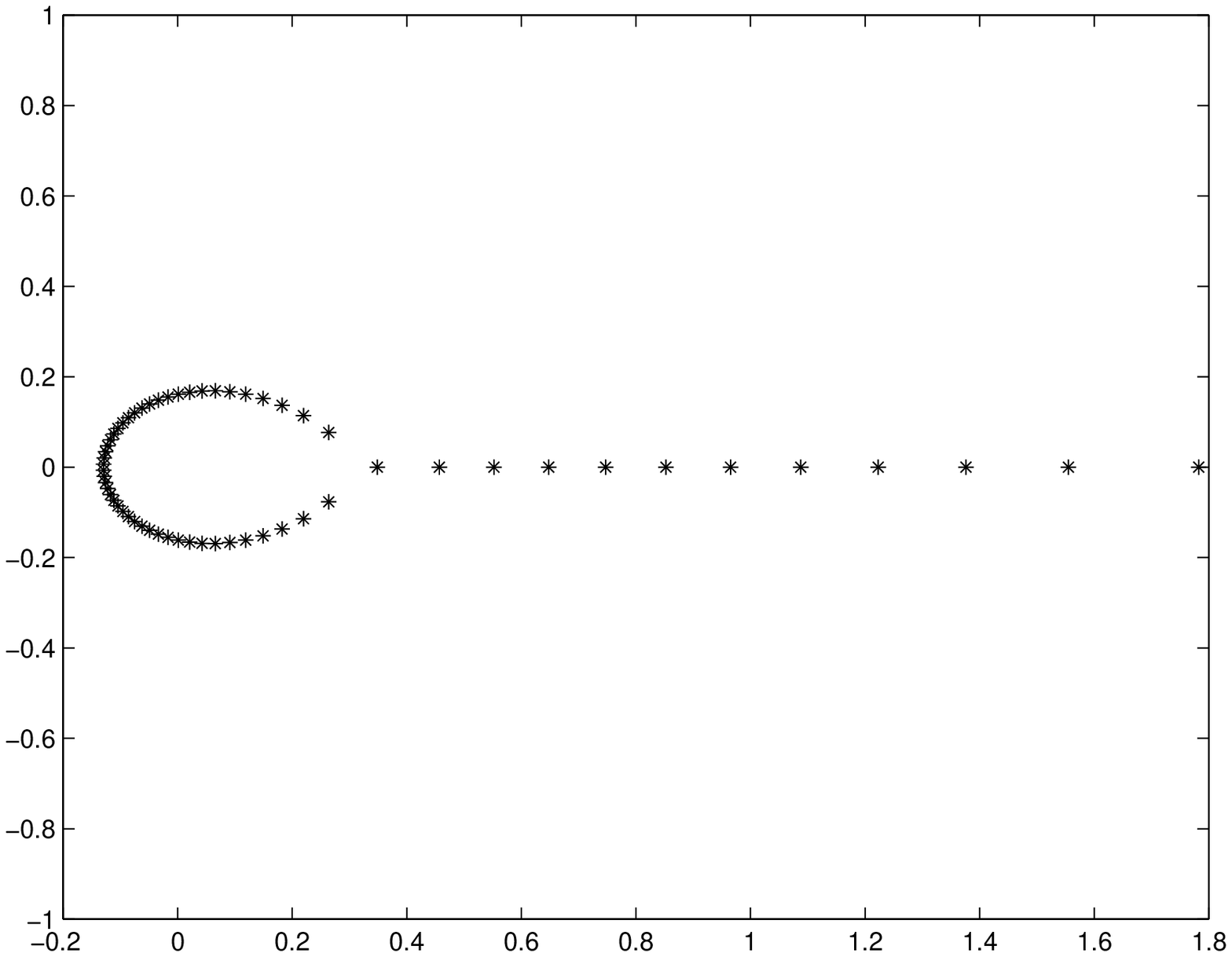}
\includegraphics[width=5cm]{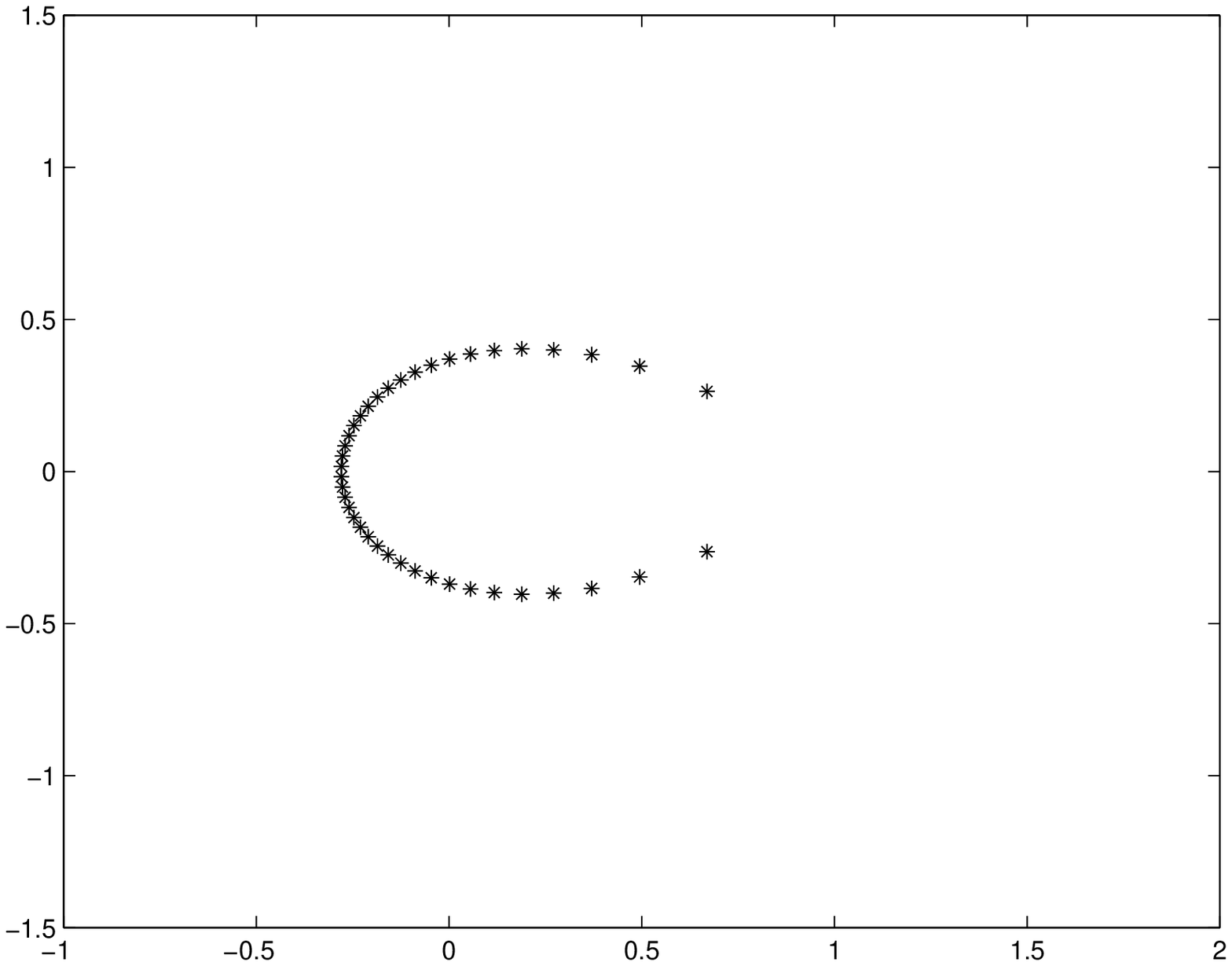}
\includegraphics[width=5cm]{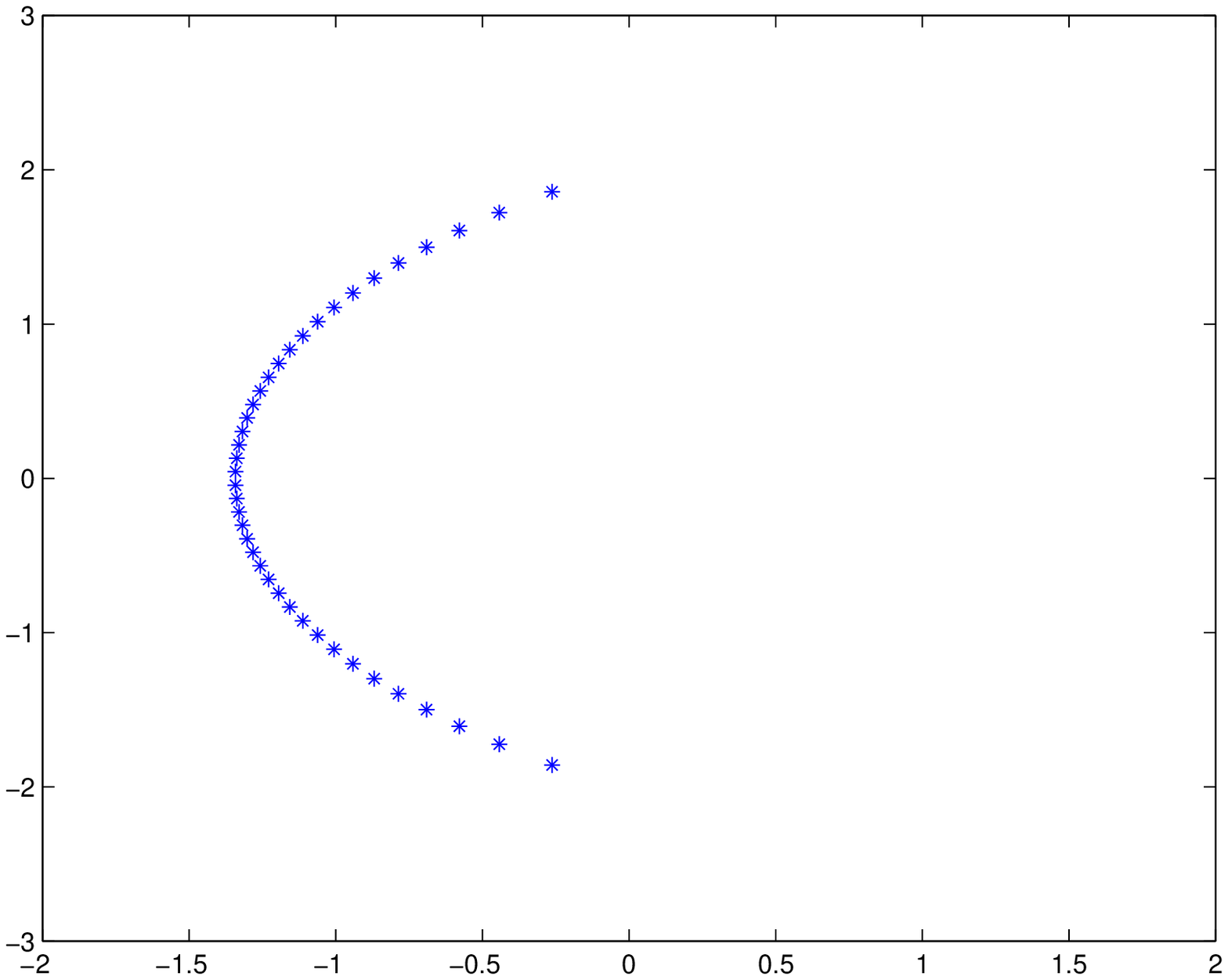}}
\caption{Zeros of Laguerre polynomials $L_{n}^{(-nA)}(nz)$ for $n = 40$
and $A = 0.81$ (left), $A = 1.01$ (middle), and $A = 2$ (right).}
\label{fig:plot1}
\end{figure}

Non-trivial limiting behavior takes place for the zeros of scaled
Laguerre polynomials $L_n^{(\alpha_n)}(nz)$ with varying negative
parameters $\alpha_n$ such that the limit
\begin{equation} \label{eq12}
     \lim_{n \to \infty} -\frac{\alpha_n}{n} =  A > 0
\end{equation}
exists.
The asymptotic location of the zeros depends on $A$
as shown in Figure \ref{fig:plot1} which is taken from \cite{KuMc},
see also \cite{MMO1}. See \cite{MMO2} for analogous plots of
zeros of Jacobi polynomials with negative parameters.

From the plots in  Figure \ref{fig:plot1} it seems clear that the
zeros accumulate along certain contours in the complex plane. For
$A > 1$, the contour is an open arc, while for $0 < A < 1$, the
contour consists of a closed loop together with an interval on the
positive real axis. In the intermediate case $A =1$, the contour
is a simple closed contour.

The case $A > 1$ is well-understood. Indeed, Saff and Varga \cite{SaVa}
showed that all of the zeros accumulate on a well-defined open arc.
Saff and Varga stated their results in terms of the zeros of the
numerator and denominator polynomials of Pad\'e approximants to the
exponential functions, but these polynomials can be identified as
Laguerre polynomials with integer parameters less than $-n$.
Uniform asymptotics for the Laguerre polynomials $L_n^{(\alpha_n)}(nz)$
in case $A > 1$ were obtained more recently
in \cite{Du}, \cite{KuMc}, and \cite{WZ}.

The case $0 < A < 1$ is less well-understood. The asymptotic
expansion of  Dunster \cite{Du} is also valid for this case, but
his results on zeros are restricted to the case $A > 1$.
Mart\'{\i}nez et al.\ \cite{MMO1} conjectured that if $(\alpha_n)$
is a sequence of non-integers such that the limit (\ref{eq12})
exists with $0 < A < 1$, then the zeros of the polynomials
$L_n^{(\alpha_n)}(nz)$ accumulate on a well-defined contour which
they identify as a trajectory of a quadratic differential. As
stated this conjecture is not true. Indeed, if every $\alpha_n$ is
an integer, then there are $|\alpha_n|$ zeros at the origin. Since
the zeros depend continuously on the parameter, it will follow
that also for $\alpha_n$ very close to the integers, some of the
zeros of the polynomials $L_n^{(\alpha_n)}(nz)$ are still close to
the origin. However, as we will show, the conjecture in
\cite{MMO1} is true provided the distance of the parameters
$\alpha_n$ to the integers is not exponentially small as $n \to \infty$.

It is the purpose of this paper to study the possible limiting behaviors
of the zeros of $L_n^{(\alpha_n)}(nz)$ as $n \to \infty$ such that (\ref{eq12})
holds with $A \in (0,1)$. We define
\begin{equation} \label{defbeta1beta2}
    \beta_1 = 2-A - 2\sqrt{1-A}, \qquad \beta_2 = 2-A + 2 \sqrt{1-A}.
\end{equation}
and
\begin{equation} \label{defRz}
    R(z) = (z-\beta_1)^{1/2} (z-\beta_2)^{1/2}, \qquad
    z \in \mathbb C \setminus [\beta_1, \beta_2],
\end{equation}
where the branches of the square roots are chosen so that $R(z) \sim z$
as $z \to \infty$.
We also put
\begin{equation} \label{defphiz}
    \phi(z) = \frac{1}{2} \int_{\beta_1}^z \frac{R(s)}{s} \; ds,
    \qquad z \in \mathbb C \setminus \left( (-\infty,0]\cup[\beta_1,\infty)\right),
\end{equation}
where the path of integration from $\beta_1$ to $z$ lies entirely in
$\mathbb C \setminus \left((-\infty,0] \cup [\beta_1,\infty)\right)$ except
for the initial point $\beta_1$.

Due to the pole of $R(s)/s$ at $s=0$ with residue $R(0) = - A$, the $\phi$
function has a jump on the negative real axis,
\begin{equation} \label{jumpphi}
    \phi_+(x) - \phi_-(x) = - A\pi i, \qquad x \in (-\infty,0),
\end{equation}
where $\phi_+$ and $\phi_-$ denote the limits from above and below, respectively.
In view of (\ref{jumpphi}), $\Re \phi$ is a well-defined harmonic function on
$\mathbb C \setminus \left(\{0\} \cup [\beta_1, \infty)\right)$.

\begin{figure}[th!]
\centerline{\includegraphics[width=8cm,height=4cm]{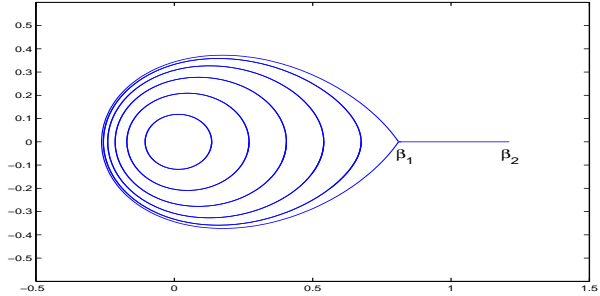}}
\caption{The interval $[\beta_1,\beta_2]$ and some of the contours $\Gamma_r$
for the value $A = 0.99$.}
\label{fig:plotGr}
\end{figure}

\begin{theorem}
For every $r \geq 0$, the following hold.
\begin{enumerate}
\item[\rm (a)] There exists a simple closed curve $\Gamma_r \subset \mathbb C \setminus
\left( \{0\} \cup (\beta_1, \infty)\right)$ encircling $0$ once in the
clockwise direction, so that
\[ \Re \phi(z) = r/2 \qquad \mbox{for } z \in \Gamma_r. \]
\item[\rm (b)]
The (a priori complex) measure
\[ d\nu_r(s) = \frac{1}{2\pi i} \frac{R(s)}{s} ds \qquad \mbox{ on } \Gamma_r \]
where $\Gamma_r$ is oriented clockwise, is a positive measure with total mass $A$.
\item[\rm (c)]
The measure
\begin{equation} \label{defmur}
    d\mu_r = d\nu_r +  \frac{\sqrt{(x-\beta_1)(\beta_2-x)}}{2\pi x}
    \chi_{[\beta_1,\beta_2]}(x) dx
\end{equation}
is a probability measure on $\Gamma_r \cup [\beta_1, \beta_2]$.
\end{enumerate}
\end{theorem}
Some of the contours $\Gamma_r$ are shown in Figure \ref{fig:plotGr}
together with the interval $[\beta_1,\beta_2]$.
\medskip

Now we consider a sequence $L_n^{(\alpha_n)}(nz)$ of rescaled Laguerre polynomials
with varying parameters $\alpha_n$ such that (\ref{eq12}) holds with $A \in (0,1)$.
We say that the zeros accumulate on a compact set $K$ if for every open set $\Omega
\supset K$, there is $n_0$ such that for every $n \geq n_0$, all zeros of
$L_n^{(\alpha_n)}(nz)$ are contained in $\Omega$. We say that a probability measure
$\mu$ is the asymptotic zero distribution if
\[ \lim_{n \to \infty} \frac{1}{n} \sum_{j=1}^n \delta_{z_{j,n}} = \mu \]
in the sense of weak$^*$ convergence of measures, see e.g.\ \cite{SaffTotik}.
Here $z_{1,n}, \ldots, z_{n,n}$
are the zeros of $L_n^{(\alpha_n)}(nz)$, counted according to multiplicity,
and $\delta_z$ is the unit Dirac measure at $z$.

\begin{theorem}
Let $(\alpha_n)$ be a sequence such that the limit {\rm (\ref{eq12})}
holds with $A \in (0,1)$. Suppose that the limit
\begin{equation} \label{defL}
    L = \lim_{n\to \infty} \left[\dist(\alpha_n, \mathbb Z)\right]^{1/n}
\end{equation}
exists.
Then the following hold.
\begin{enumerate}
\item[\rm (a)] If $L=0$, then the zeros of $L_n^{(\alpha_n)}(nz)$
accumulate on $\{ 0 \} \cup [\beta_1, \beta_2]$ and
\begin{equation} \label{MarPas}
    A \delta_0 + \frac{\sqrt{(x-\beta_1)(\beta_2-x)}}{2\pi x} \chi_{[\beta_1,\beta_2]}(x) dx
\end{equation}
is the asymptotic zero distribution.
\item[\rm (b)] If $L = e^{-r}$ with $0 \leq r < \infty$, then the zeros
of $L_n^{(\alpha_n)}(nz)$ accumulate on $\Gamma_r \cup [\beta_1, \beta_2]$
and the measure $\mu_r$ given in {\rm(\ref{defmur})} is the asymptotic zero distribution.
\end{enumerate}
\end{theorem}

\begin{remark}
The case $L=1$ in Theorem 1.2 represents the typical case, in
the sense that if a sequence $(\alpha_n)$ would be chosen randomly
according to some reasonable distribution, then then
\[ \lim_{n\to\infty} \left[\dist(\alpha_n,\mathbb Z)\right]^{1/n} = 1 \]
with probability one. So in the typical case, the zeros cluster
on $[\beta_1, \beta_2]$ and along the outer curve $\Gamma_0$.
The case of limit $L < 1$ in (\ref{defL}) is more special, since
for that to happen the parameters $\alpha_n$ should be very close
to integers.

To illustrate this, we have plotted in Figure \ref{fig:plot1b} the
zeros of $L_{n}^{(-nA)}$ with $n = 40$, for the values $A = 0.81$ and $A = 0.799999975$.
Then the distance of $\alpha_n = -nA$ to the integers is $0.4$ and
$10^{-6}$, respectively. There is a clear separation in the plot on
the right between the zeros on the interval $[\beta_1,\beta_2]$ and
the cluster of zeros around the origin. For the plot on the right,
$\beta_1 = 0.31\cdots$ and $\beta_2 = 2.09\cdots$.
The rightmost zero in the cluster has real part $0.14\cdots$.
\begin{figure}[th!]
\centerline{\includegraphics[width=5cm]{plot1d.eps} \qquad
\includegraphics[width=5cm]{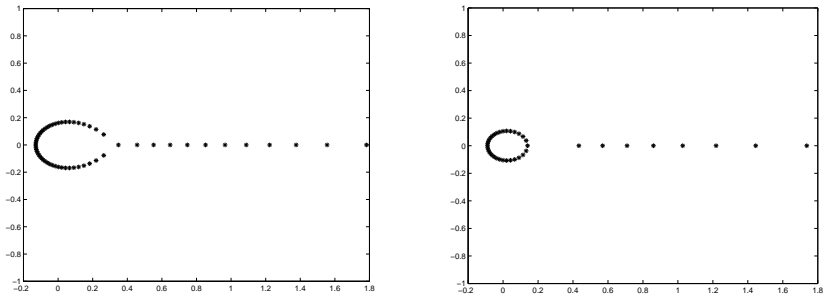}}
\caption{Zeros of Laguerre polynomials $L_{n}^{(-nA)}(nz)$ for $n = 40$
and $A = 0.81$ (left), $A = 0.799999975$ (right).
The two cases correspond to $\dist(\alpha_n,\mathbb Z) = 0.4$ and
$\dist(\alpha_n,\mathbb Z) = 10^{-6}$, respectively.}
\label{fig:plot1b}
\end{figure}
\end{remark}

\begin{remark}
Theorem 1.2 also applies to subsequences. So, if $(\alpha_n)$ is a
sequence such that the limit (\ref{defL}) does not exist, then we
can pass to a subsequence such that (\ref{defL}) exists
along that subsequence. Then we find an asymptotic distribution
of zeros for any subsequence such that the limit (\ref{defL}) exists.
By constructing special sequences we can arrange it that different
subsequences have different asymptotic zero distributions, and it
is even possible to construct examples in which every measure
$\mu_r$ appears as the asymptotic zero distribution along some
subsequence.
\end{remark}

\begin{remark}
For fixed $n$ and decreasing  $\alpha < 0$, the dynamics of the
zeros of $L_n^{(\alpha)}$  is similar to the dynamics of the zeros
of certain hypergeometric polynomials  studied recently by Driver
and Duren \cite{DD1,DD2}. For the Laguerre polynomials, the zeros
are all real and positive if $\alpha > -1$. For $\alpha=-1$, there
is one zero at the origin. This zero becomes negative if $\alpha$
becomes smaller than $-1$. As $\alpha$ decreases, the smallest
zero reaches a minimum and then starts to increase again. For
$\alpha = -2$, it is at the origin again, and simultaneously the
second smallest zero has reached the origin. A collision takes
place, and for $\alpha$ somewhat smaller than $-2$, the two zeros
separate and leave the real axis, moving into the complex plane.
When $\alpha$ has reached the next integer value $-3$ they both
return to the origin, and a third zero coming from the positive
real axis joins them in a three-fold collision. As $\alpha$
decreases from $-3$ to $-4$, the three zeros move into the complex
plane, one on the negative real axis, and two away from the axis,
coming again to the origin for a new collision at $\alpha = -4$,
where they are joined by another zero coming from the positive
real axis. And so on, until all $n$ zeros have left the positive
real axis. At $\alpha = -n$, all $n$ zeros meet at the origin. For
$\alpha < -n$, they separate again, and they all move into the
complex plane. If $n$ is odd, then one zero remains on the
negative real axis.
\end{remark}

\begin{remark}
In statistical literature, see \cite{Bai,Dette} and the references
cited therein, the measure (\ref{MarPas}) is
known as the Marchenko-Pastur distribution \cite{MarPas}.
It appears as the
limiting spectral distribution of Wishart matrices, which are
random covariance matrices $\frac{1}{p} X X^T$ where
$X$ is an $n\times p$ matrix with independent, identically distributed
entries. The measure (\ref{MarPas}) corresponds to the limit
$n,p \to \infty$, $n/p \to 1-A \in (0,1)$.

In the theory of free probability, the measure (\ref{MarPas})
is also known as the free Poisson distribution, see \cite{HiaiPetz}
and the references cited therein.
\end{remark}

We prove Theorem 1.1 in section 2. The rest of the paper is devoted
to the proof of Theorem 1.2. It is based on a characterization of Laguerre
polynomials by a Riemann-Hilbert problem \cite{FIK,KuMc},
and an asymptotic analysis of the Riemann-Hilbert problem with
the steepest descent / stationary phase method introduced by
Deift and Zhou \cite{DZ}, and developed further in
\cite{BleherIts,DKMVZ1,DKMVZ2,DVZ,KrMc,KMVV}. See \cite{Deift} for an
introduction to this method, with a special emphasis on the
applications to orthogonal polynomials and random matrices.
In \cite{KuMc} we used the Riemann-Hilbert method to study
Laguerre polynomials with varying negative parameters $\alpha_n$
such that $\lim_{n} \alpha_n/n  < -1$, which is the case $A > 1$.
In this case, the zeros accumulate along an open contour as shown
in Figure 1. The determination of this contour is a first step
in the analysis in \cite{KuMc}.  Riemann-Hilbert analysis in
which a contour selection method was required appeared also in the
recent papers \cite{Aptekarev, BDMMZ,KMM}.

In the Riemann-Hilbert method, the endpoints of the contour need
special attention. Following \cite{DKMVZ1,DKMVZ2}, the local
analysis around the endpoints is done in \cite{KuMc} with the help of
Airy functions by constructing a so-called Airy parametrix.
The analysis for the case $0 < A < 1$ proceeds along similar lines.
However, the fact that the zeros accumulate along a more complicated
contour causes some new features. In particular we will be dealing
with a closed loop $\Gamma_0$ together with an interval
$[\beta_1,\beta_2]$ on the real line. The point $\beta_1$ where the
closed loop and the interval intersect is a new special point.
The surprising fact, however, is that after suitable transformations,
the transformed Riemann-Hilbert problem is very similar
to the transformed Riemann-Hilbert considered
in \cite{DKMVZ2}. It then follows that the local analysis around
$\beta_1$ can be done with the Airy parametrix as well.

The Riemann--Hilbert method gives strong and uniform asymptotics
of the Laguerre polynomials $L_n^{(\alpha_n)}(nz)$ in the
whole complex plane. This is really much more than what is needed
for the proof of Theorem 1.2. Since we do not have an application
for strong asymptotics (although we could envision there is one),
we have not pursued this matter any further. As a sample of the
kind of results that could be obtained, we state here the
asymptotic result in the oscillatory region $(\beta_1, \beta_2)$
on the real line. Suppose that the limit (\ref{eq12}) holds
with $A \in (0,1)$. For every $n \in \mathbb N$, we put
$A_n = - \frac{\alpha_n}{n}$, so that $A_n \to A$ as $n \to \infty$.
Let $\beta_{1,n}$ and $\beta_{2,n}$ be the values (\ref{defbeta1beta2})
corresponding to $A_n$ instead of $A$. Then we have, uniformly
for $x$ in compact subsets of $(\beta_1, \beta_2)$,
\begin{eqnarray} \nonumber
  \lefteqn{  L_n^{(\alpha_n)}(nx) = \frac{(-1)^n}{n!}
        \frac{2}{\sqrt{\beta_{2,n}-\beta_{1,n}}}
        \frac{1}{(\beta_{2,n}-x)^{1/4}(x-\beta_{1,n})^{1/4}}} \\
        & & \nonumber
        \left\{ \cos \left(n\pi \int_{\beta_{2,n}}^x
            \frac{\sqrt{(s-\beta_{1,n})(\beta_{2,n}-s)}}{2\pi s} ds
            + \frac{1}{2} \arcsin \frac{2z-(\beta_{2,n} + \beta_{1,n})}{\beta_{2,n}-\beta_{1,n}}
            \right) \left(1+ O\left(\frac{1}{n}\right)\right) \right. \\
    & &  \label{oscilasymp} \left.
        + \sin\left(n\pi \int_{\beta_{2,n}}^x
            \frac{\sqrt{(s-\beta_{1,n})(\beta_{2,n}-s)}}{2\pi s} ds
            - \frac{1}{2} \arcsin \frac{2z-(\beta_{2,n} + \beta_{1,n})}{\beta_{2,n}-\beta_{1,n}}
            \right) O\left(\frac{1}{n}\right)
    \right\}
\end{eqnarray}

\section{Proof of Theorem 1.1}
\setcounter{equation}{0}

\begin{proof}
The curves where $\Re \phi(z)$ is constant are known as trajectories
of the quadratic differential
\begin{equation} \label{defquaddiff}
    -\frac{R(z)^2}{z^2} dz^2 = -\frac{(z-\beta_1)(z-\beta_2)}{z^2} dz^2,
\end{equation}
which has simple zeros at $\beta_1$ and $\beta_2$ and a double pole at $0$,
see \cite{Strebel}.
From the local structure at a simple zero of the
quadratic differential, see \cite[\S 7]{Strebel}, we know that
three trajectories emanate from $\beta_1$ at equal angles $2 \pi/3$.
One of these is the interval $[\beta_1, \beta_2]$. Let $\Gamma_0$ be the trajectory emanating
from $\beta_1$ at the angle $-2\pi/3$, and let $z(t)$ be the arclength
parametrization of $\Gamma_0$ with $z(0) = \beta_1$.
So $\Re[ \phi(z(t))] = 0$ for every $t$, which implies in particular
that $\Gamma_0$ does not approach the origin.
Since $\Gamma_0$ leaves $\beta_1$ at an angle $-2\pi/3$, we
have for small $t >0$ (and  thus for
$z(t)$ close to $\beta_1$), that $|z(t)|$ decreases as $t$ increases.
Then $|z(t)|$ continues to decrease until for some $t_0$, we have
\begin{equation} \label{criticaltraj1}
    \Re \left[ \frac{z'(t_0)}{z(t_0)} \right ] = 0.
\end{equation}
Since $\Re[\phi(z(t))] = 0$ for every $t$, we find by differentiating
(\ref{defphiz}),
\begin{equation} \label{criticaltraj2}
    \Re \left[ z'(t) \frac{R(z(t))}{z(t)} \right] = 0
\end{equation}
for every $t$. Thus for $t = t_0$ we have from (\ref{criticaltraj1}) and
(\ref{criticaltraj2}) that $R(z(t_0)) \in \mathbb R$.
Since $R(z)$ is only real for $z$ real, it follows that $z(t_0)$ is real,
So $|z(t)|$ decreases until $\Gamma_0$ hits the real axis.
Since $z(0) = \beta_1$, we have $|z(t_0)| < \beta_1$.
From (\ref{defphiz}) it is clear that $\phi(x)$ is real and strictly
positive for $x \in (0, \beta_1)$, since $R(s) < 0$ for $s \in (-\infty,\beta_1)$.
Thus $z(t_0) < 0$.
So we have a curve in the lower half-plane from $\beta_1$ to $z(t_0)$
where $\Re \phi = 0$. By symmetry, we also have $\Re \phi = 0$ on the
mirror image of this curve. Then the curve and its mirror image form
the full trajectory $\Gamma_0$ which encircles the origin. This proves
part (a) of Theorem 1.1 for $r = 0$.

See Figure \ref{fig:traject} for a plot of $\Gamma_0$ for two values
of $A$. Also plotted are the other contours where $\Re \phi(z) = 0$ (solid
lines). The dotted lines are the curves where either  $\phi(z)$ or
\begin{equation} \label{deftildephiz}
    \tilde{\phi}(z) = \frac{1}{2} \int_{\beta_2}^z \frac{R(s)}{s} \, ds
\end{equation}
are real.
\begin{figure}[th!]
\centerline{\includegraphics[width=7cm]{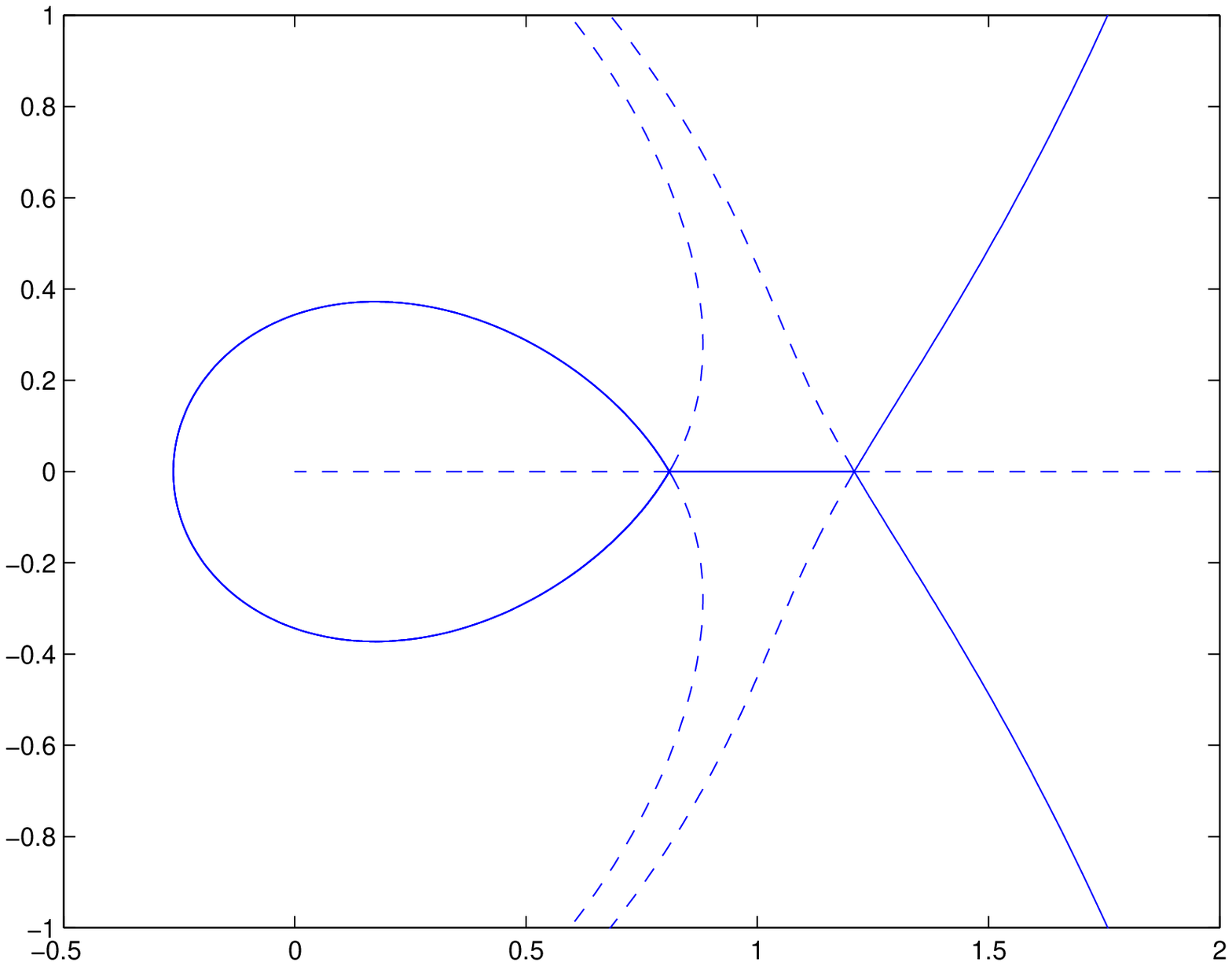}
    \includegraphics[width=7cm]{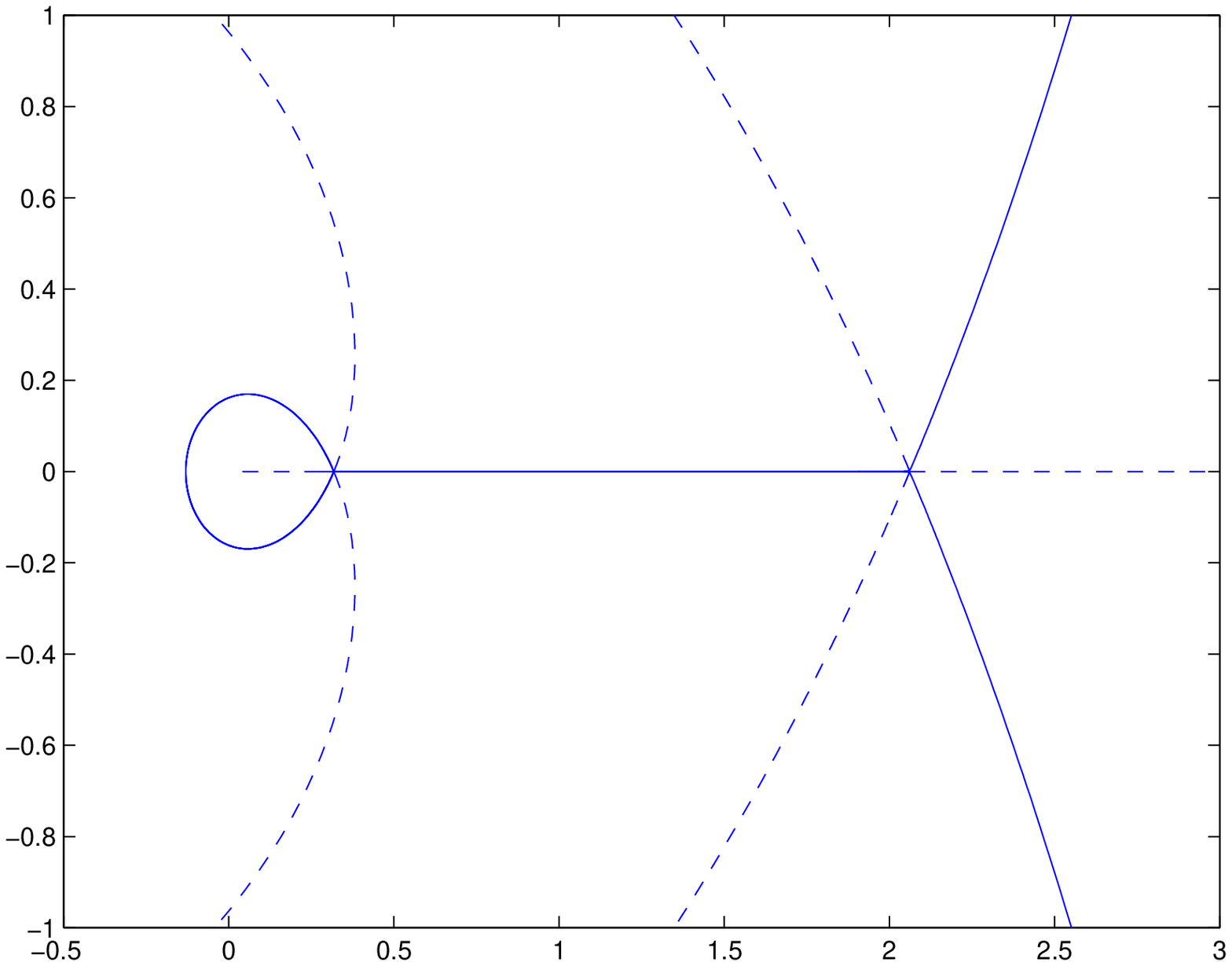}}
\caption{Curves where $\phi$ is purely imaginary
(solid lines), or where either $\phi$ or $\tilde{\phi}$ are real (dotted lines).
The closed  loop is $\Gamma_0$.
The values for $A$ are $0.99$ (right) and $0.81$ (left).}
\label{fig:traject}
\end{figure}
\medskip

Let $\Omega_0$ be the region that is encircled by $\Gamma_0$. Then
$\Re \phi$ is harmonic on $\Omega_0 \setminus \{ 0 \}$ with zero boundary
values on $\Gamma_0$. From (\ref{defphiz}) it follows that $\Re \phi$
behaves like $\frac{A}{2} \log |z|$ as $z \to 0$.
Thus $(2/A) \Re \phi$ is the Green function of $\Omega_0$ with pole at $0$.
Since $\Gamma_0$ is a simple closed curve, it then follows that for
every $r >0$ the level curve
\[ \Gamma_r = \{ z \in \Omega_0 \mid \Re \phi(z) = r/2 \} \]
is a simple closed curve encircling the origin. This proves part (a)
for every $r \geq 0$.
\medskip

For part (b) we first note that by Cauchy's formula
\[ \int_{\Gamma_r} d\nu_r(s) = \frac{1}{2\pi i} \int_{\Gamma_r} \frac{R(s)}{s} ds
= - R(0). \]
We have a $-$sign since $\Gamma_r$ is oriented in clockwise direction.
Since $R(0) = -A$, we get $\int_{\Gamma_r} d\nu_r = A > 0$.
If $z(t)$ is arclength parametrization of $\Gamma_r$ with $z(0) = \beta_1$, then
\[ t \mapsto \int_{\beta_1}^{z(t)} d\nu_r = \frac{1}{2\pi i} \int_{\beta_1}^{z(t)} \frac{R(s)}{s} ds \]
is differentiable with derivative $\frac{1}{2\pi i} \frac{R(z(t))}{z(t)} z'(t)$,
which is real by the definition of $\Gamma_r$, and non-zero since $R$ is non-zero
on $\Gamma_r$. So $\int_{\beta_1}^{z(t)} d\nu_r$ can only increase along $\Gamma_r$
and it follows that $\nu_r$ is a positive measure with total mass $A$. This proves part (b).
\medskip

Then for part (c) it is enough to observe that
\[ \int_{\beta_1}^{\beta_2} \frac{\sqrt{(\beta_1-x)(x-\beta_2)}}{x} \; dx = 1-A, \]
which can be shown by direct calculation, or by a contour deformation argument.
\end{proof}

\section{Proof of Theorem 1.2}
\setcounter{equation}{0}
\subsection{Preliminaries}
We assume that the sequence $(\alpha_n)$ is as in Theorem 1.2. We put
\begin{equation} \label{defAn}
    A_n = - \frac{\alpha_n}{n}.
\end{equation}
so that $\lim\limits_{n \to \infty} A_n = A$. Since $A \in (0,1)$, we may assume
without loss of generality that $A_n \in (0,1)$ for all $n$, so that
$\alpha_n \in (-n,0)$.
We may also assume that $\alpha_n \not\in \mathbb Z$. Indeed,
if $\alpha_n \in \mathbb Z$, then by (\ref{eq11a}) the zeros of
$L_n^{(\alpha_n)}(nz)$ are the zeros of $L_{n+\alpha_n}^{(-\alpha_n)}(nz)$,
together with a zero at $0$ of multiplicity $|\alpha_n|$.
From the asymptotic behavior of zeros of Laguerre polynomials
with varying positive parameters, see \cite{DS,Gaw1,Gaw2,KuVA},
part (a) of Theorem 1.2 follows, in case the $\alpha_n$ are integers.
So we may and do assume that $\alpha_n \not\in \mathbb Z$.

For every $n$, we use a subscript $n$ to denote the notions (\ref{defbeta1beta2}),
(\ref{defRz}), (\ref{defphiz}), but with $A$ replaced by $A_n$. Thus
\begin{equation} \label{defbeta1nbeta2n}
    \beta_{1,n} = 2-A_n - 2\sqrt{1-A_n}, \qquad \beta_{2,n} = 2-A_n + 2 \sqrt{1-A_n}.
\end{equation}
\begin{equation} \label{defRnz}
    R_n(z) = (z-\beta_{1,n})^{1/2} (z-\beta_{2,n})^{1/2}, \qquad
    z \in \mathbb C \setminus [\beta_{1,n}, \beta_{2,n}],
\end{equation}
and
\begin{equation} \label{defphinz}
    \phi_n(z) = \frac{1}{2} \int_{\beta_{1,n}}^z \frac{R_n(s)}{s} \; ds,
    \qquad z \in \mathbb C \setminus \left( (-\infty,0]\cup[\beta_{1,n},\infty)\right).
\end{equation}
We also put
\begin{equation} \label{deftildephinz}
    \tilde{\phi}_n(z) = \frac{1}{2} \int_{\beta_{2,n}}^z \frac{R_n(s)}{s} \; ds,
    \qquad z \in \mathbb C \setminus (-\infty, \beta_{2,n}],
\end{equation}
where the path of integration from $\beta_{2,n}$ to $z$ lies entirely in
$\mathbb C \setminus (-\infty,\beta_{2,n}]$, except for the initial point $\beta_{2,n}$.
We denote by $\Gamma_0^{(n)}$ the simple closed contour around $0$ given by
$\Re \phi_n(z) = 0$, whose existence is guaranteed by Theorem 1.1(a).

\subsection{The Riemann-Hilbert problem}

For every unbounded contour $\Sigma$ in $\mathbb C \setminus [0,\infty)$
from a point $+\infty -iy$ to $+\infty + iy$ for some $y > 0$,
we have the orthogonality relation
\begin{equation} \label{eq22}
    \int_{\Sigma} L_n^{(- n A_n)}(z) z^k z^{-n A_n} e^{-z} dz = 0,
    \quad \mbox{ for } k = 0, 1, \ldots, n-1,
\end{equation}
where the branch of $z^{-n A_n}$ is taken with a cut along
the positive real axis. In addition, since $n A_n \not\in \mathbb Z$,
\begin{equation} \label{eq23}
    \int_{\Sigma} L_n^{(-n A_n)}(z) z^n z^{-nA_n} e^{-z} dz \neq 0.
\end{equation}
The relations (\ref{eq22}) and (\ref{eq23}) were proved in
\cite[Lemma 2.1]{KuMc}. They follow from the Rodrigues formula for
Laguerre polynomials.
\begin{figure}[th!]
\centerline{\includegraphics[width=8cm]{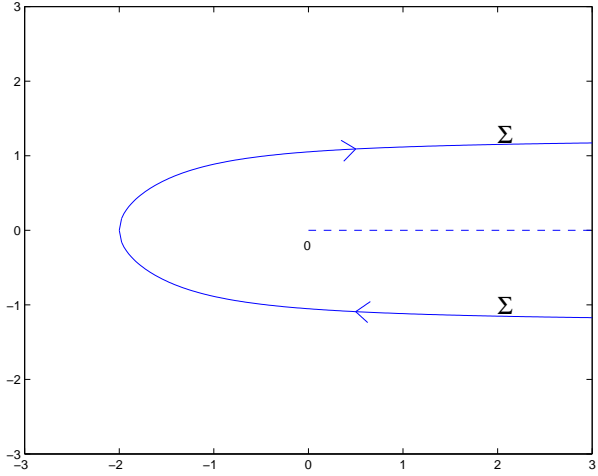}}
\caption{Example of contour $\Sigma$}
\label{fig:sigmaorig}
\end{figure}

The rescaled monic polynomial
\begin{equation} \label{eq24}
    P_n(z) = \frac{n!}{(-n)^n} L_n^{(-n A_n)}(nz)
\end{equation}
is then characterized by a Riemann-Hilbert problem
on a contour $\Sigma$ as shown in Figure~\ref{fig:sigmaorig}, oriented
from $+\infty-iy$ to $+\infty+iy$. For any oriented contour $\Sigma$, we follow
the usual convention that the $+$-side is on the left and the $-$-side is
on the right while traversing the contour. We say that a function $Y$ defined
on $\mathbb C \setminus \Sigma$ has boundary values on $\Sigma$ if  for every
$z \in \Sigma$, the two limits
\[ Y_+(z) = \lim_{{\rm+side} \ni z'\to z} Y(z'), \qquad
   Y_-(z) = \lim_{{\rm-side} \ni z' \to z} Y(z')
\]
exist.

\subsubsection*{Riemann-Hilbert problem for $Y$:}
Let $\Sigma$ be as above.
The problem is to determine a $2 \times 2$
matrix valued function $Y : \mathbb C \setminus \Sigma \to \mathbb C^{2 \times 2}$
such that the following hold.
\begin{enumerate}
\item[(a)] $Y(z)$ is analytic for $z \in \mathbb C \setminus \Sigma$,
\item[(b)] $Y(z)$ possesses continuous boundary values for $z \in \Sigma$,
denoted by $Y_+(z)$ and $Y_-(z)$,  and
\[ Y_+(z) = Y_-(z)
    \left( \begin{array}{cc} 1 & z^{-n A_n} e^{-nz} \\ 0 & 1 \end{array} \right),
    \qquad \mbox{ for } z \in \Sigma,
\]
\item[(c)] $Y(z)$ has the following behavior as $z \to \infty$:
\[  Y(z) = \left( I + O\left(\frac{1}{z}\right) \right)
    \left( \begin{array}{cc} z^n & 0 \\ 0 & z^{-n} \end{array} \right)
    \qquad \mbox{ as } z \to \infty, \ z \in \mathbb C \setminus \Sigma.
\]
\end{enumerate}
The Riemann-Hilbert problem has a unique solution,
see \cite[Proposition 2.2]{KuMc}, given by
\begin{equation} \label{solutionY}   Y(z) = \left( \begin{array}{ccc}
    P_n(z) & \frac{1}{2\pi i}  \int\limits_{\Sigma}
    \frac{P_n (\zeta) \zeta^{-nA_n} e^{-n\zeta}}{\zeta - z}
     d\zeta \\[10pt]
    Q_{n-1}(z) & \frac{1}{2\pi i}  \int\limits_{\Sigma}
     \frac{Q_{n-1}(\zeta) \zeta^{-nA_n} e^{-n\zeta}}{\zeta -z } d\zeta
     \end{array} \right)
\end{equation}
where $Q_{n-1}(z)= d_n L_{n-1}^{(-n A_n)}(nz)$ is a polynomial of degree $n-1$,
and the constant $d_n$ is chosen such that $Y_{22}(z)$ behaves exactly like $z^{-n}$ as
$z \to \infty$.
\medskip

\begin{figure}[th!]
\centerline{\includegraphics[width=8cm]{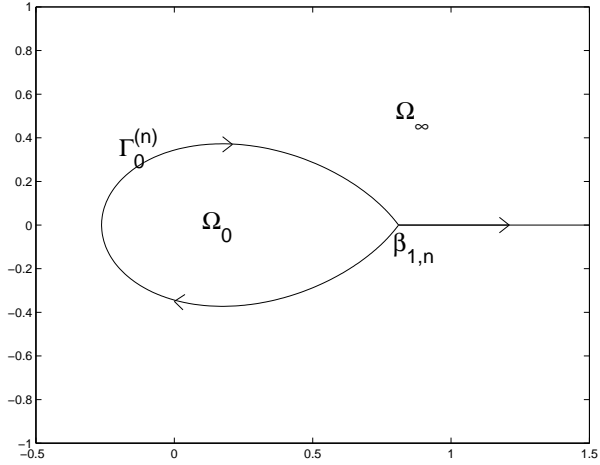}}
\caption{Contour $\Sigma^{(n)} = \Gamma_0^{(n)} \cup [\beta_{1,n}, \infty)$}
\label{fig:sigma}
\end{figure}

Anticipating the distribution of zeros as shown in Figure \ref{fig:plot1} for
the case $0 < A < 1$, we modify the Riemann-Hilbert problem for $Y$
as follows. We let part of the contour $\Sigma$ come to the positive real
axis and put
\[ \Sigma = \Sigma^{(n)} = \Gamma_0^{(n)} \cup [\beta_{1,n}, \infty). \]
The curve $\Gamma_0^{(n)}$ is oriented in the clockwise direction,
and $[\beta_{1,n},\infty)$ is oriented from left to right. The unbounded component
of $\mathbb C \setminus \Sigma^{(n)}$ is denoted by $\Omega_{\infty}$ and the bounded
component by $\Omega_0$ as indicated in Figure \ref{fig:sigma}.

Then it immediately follows that $P_n$ is characterized by the following
Riemann-Hilbert problem.
\subsubsection*{Modified Riemann-Hilbert problem for $Y$:}
Let $\Sigma^{(n)}$ be as described above.
The problem is to determine a $2 \times 2$
matrix valued function $Y : \mathbb C \setminus \Sigma \to \mathbb C^{2 \times 2}$
such that the following hold.
\begin{enumerate}
\item[(a)] $Y(z)$ is analytic for $z \in \mathbb C \setminus \Sigma^{(n)}$,
\item[(b)] $Y(z)$ possesses continuous boundary values for $z \in \Sigma^{(n)}$,
denoted by $Y_+(z)$ and $Y_-(z)$,  and $Y_+(z) = Y_-(z) V_Y(z)$ where
\begin{align}
    \label{jumpVY1}
    V_Y(z) &= \left( \begin{array}{cc} 1 & z^{-n A_n} e^{-nz} \\ 0 & 1 \end{array} \right),
    && \mbox{for } z \in \Gamma_0^{(n)}, \\[10pt]
    \label{jumpVY2}
    V_Y(z) &=
    \left(\begin{array}{cc} 1 & c_n z^{-n A_n} e^{-nz} e^{-nA_n\pi i} \\
    0 & 1 \end{array} \right),
    && \mbox{for } z \in (\beta_{1,n},\infty),
\end{align}
and $c_n = 2i \sin(nA_n\pi)$.
\item[(c)] $ \ds Y(z) = \left( I + O\left(\frac{1}{z}\right) \right)
    \left( \begin{array}{cc} z^n & 0 \\ 0 & z^{-n} \end{array} \right)$
     as $ z \to \infty$.
\end{enumerate}
Then as before we have that the modified Riemann-Hilbert problem
for $Y$ has a unique solution and that
\begin{equation} \label{Y11isPn}
    Y_{11}(z) = P_n(z).
\end{equation}

The constant
\begin{equation}\label{defcn}
    c_n = 2i \sin(nA_n\pi)
\end{equation}
appearing in part (b) will continue to play a role later on.
Since $\alpha_n$ is not an integer, $c_n$
is different from zero, and it can be used as a measure for the distance
of $\alpha_n$ to the integers. In fact, we have that
\begin{equation} \label{Lincn}
    L = \lim_{n \to \infty} |c_n|^{1/n},
\end{equation}
where $L$ is given by (\ref{defL}).

\subsection{First transformation $Y \mapsto U$}

The first transformation of the Riemann-Hilbert problem is
based on the construction of a so-called $g$-function.
\begin{definition}
We define the function $g_n$ by
\begin{equation} \label{defgnz}
    g_n(z) = \frac{1}{2\pi i} \int_{\Gamma_0^{(n)} \cup [\beta_{1,n}, \beta_{2,n}]}
    \log(z-s) \frac{R_{n,+}(s)}{s} \, ds,
    \qquad z \in \mathbb C \setminus \Sigma^{(n)}.
\end{equation}
For $s \in [\beta_{1,n},\beta_{2,n}]$, the logarithm $\log(z-s)$ is taken with
a cut along $[s,\infty)$, for $s \in \Gamma_0^{(n)}$ the branch of the logarithm
$\log(z-s)$ is chosen with a cut along the part of $\Gamma_0^{(n)}$ that starts at $s$
continues along $\Gamma_0^{(n)}$ to the right until $\beta_{1,n}$, and then consists
of $[\beta_{1,n}, \infty)$.
\end{definition}

Thus $g_n$ is the log-transform of the probability measure $\mu_0$, see (\ref{defmur}),
but corresponding to the value $A_n$ instead of $A$. Note that $g_n$ is defined and
analytic in $\mathbb C \setminus \Sigma^{(n)}$.

The following proposition describes the jump properties of $g_n$. They follow
by taking into account the various branches of the logarithm, see also \cite{KuMc}.
We use $\mathbb C^+$ and $\mathbb C^-$ to denote the upper and lower half-plane,
respectively.
\begin{proposition}
\begin{enumerate}
\item[\rm (a)]  There is a constant $\ell_n$ such that
\begin{equation} \label{eq46}
    g_n(z) = \frac{1}{2} \left( A_n \log z + z + \ell_n \right) +
    \left\{ \begin{array}{ll}
    \mp \frac{1}{2} A_n\pi i  - \phi_n(z), & \quad \mbox{ for } z \in \Omega_{\infty} \cap \mathbb C^{\pm}, \\[10pt]
    \pm \frac{1}{2} A_n\pi i  + \phi_n(z), & \quad \mbox{ for } z \in \Omega_0 \cap \mathbb C^{\pm}.
    \end{array} \right.
\end{equation}
Here $\log z = \log|z| + i \arg z$  with $\arg z \in [0,2\pi)$. The number $\ell_n$ is
explicitly given by $\ell_n = 2g_n(x_n)-(A_n\log x_n + x_n)$, where $x_n$
is the intersection of $\Gamma_0^{(n)}$ with the negative real axis.
\item[\rm (b)]
We have the following jump for $g_n$ on $\Sigma^{(n)}$,
\begin{equation} \label{eq47}
g_{n,+}(z) - g_{n,-}(z) =
    \left\{ \begin{array}{ll}
    -2\phi_n(z) \mp A_n \pi i & \quad \mbox{for } z \in \Gamma_0^{(n)} \cap \mathbb C^{\pm}, \\[10pt]
    - 2\phi_{n,+}(z) - 2A_n \pi i \, = \, 2\phi_{n,-}(z) - 2A_n \pi i
    & \quad \mbox{for } z \in (\beta_{1,n},\beta_{2,n}], \\[10pt]
    -2 \pi i & \quad \mbox{for } z \in [\beta_{2,n},\infty).
    \end{array} \right.
\end{equation}
\item[\rm (c)] We have, with the same constant $\ell_n$ as in part {\rm (a)},
\begin{equation} \label{eq48}
    g_{n,+}(z) + g_{n,-}(z) =
    \left\{ \begin{array}{ll}
    A_n \log z + z + \ell_n, & \qquad \mbox{for } z \in \Gamma_0^{(n)}, \\[10pt]
    A_n \log z + z + A_n \pi i + \ell_n, & \qquad \mbox{for } z \in [\beta_{1,n},\beta_{2,n}], \\[10pt]
    A_n \log z + z + A_n \pi i + \ell_n - 2\tilde{\phi}_n(z), &
    \qquad \mbox{for } z \in (\beta_{2,n},\infty).
    \end{array} \right.
\end{equation}
\end{enumerate}
\end{proposition}
\begin{proof}
Taking the derivative of (\ref{defgnz}) we get that
\[ g_n'(z) = \frac{1}{2\pi i} \int_{\Gamma_0^{(n)}} \frac{1}{z-s} \frac{R_n(s)}{s} ds
    + \frac{1}{2\pi i} \int_{[\beta_1,\beta_2]} \frac{1}{z-s} \frac{R_{n,+}(s)}{s} ds. \]
The two integrals can be calculated using contour deformation and Cauchy's formula.
The result is that
\begin{equation} \label{eq49}
    g_n'(z) = \left\{ \begin{array}{ll}
    \frac{1}{2} \left( -\frac{R_n(z)}{z} + \frac{A_n}{z} + 1 \right), & \qquad
    \mbox{for } z \in \Omega_{\infty}, \\[10pt]
    \frac{1}{2} \left( \frac{R_n(z)}{z} + \frac{A_n}{z} + 1 \right), & \qquad
    \mbox{for } z \in \Omega_0.
    \end{array} \right.
\end{equation}
We integrate $g_n'$ from $x_n$ to $z$ along a path in
$\mathbb C \setminus (\Sigma^{(n)} \cup \mathbb R)$,
except for the initial point $x_n$.
We obtain from (\ref{eq49}) and the definition of $\ell_n$,
\begin{eqnarray*}
     g_n(z) & = & g_n(x_n) + \int_{x_n}^z g_n'(s) ds \\
     & = & \frac{1}{2} (A_n \log z + z + \ell_n)
        \left\{ \begin{array}{ll}
        - \frac{1}{2} \int_{x_n}^z \frac{R_{n}(s)}{s} ds & \mbox{ for } z \in \Omega_{\infty}, \\[10pt]
        + \frac{1}{2} \int_{x_n}^z \frac{R_n(s)}{s} ds & \mbox{ for } z \in \Omega_0,
        \end{array} \right.
\end{eqnarray*}
where we also used the fact that $g_{n,+}(x_n) = g_{n,-}(x_n)$
because of  the choice of the branches of the logarithm
in (\ref{defgnz}). [Here $g_{n,+}(x_n)$ and $g_{n,-}(x_n)$ denote the limiting values
of $g_n(z')$ as $z' \to x_n$ from within $\Omega_{\infty}$ and $\Omega_{0}$, respectively.]
By the definition (\ref{defphinz}) of $\phi_n(z)$, we have
\[ \frac{1}{2} \int_{x_n}^z \frac{R_n(s)}{s} ds =
    \phi(z) - \phi_{\pm}(x_n) \qquad \mbox{ for } z \in \mathbb C^{\pm}, \]
and (\ref{eq46}) follows, since $\phi_{\pm}(x_n) = \mp \frac{1}{2} A_n \pi i$.
This proves part (a) of the proposition. Parts (b) and (c) follow immediately
from part (a).
\end{proof}

With the $g$-function, the constant $\ell_n$, and the number $c_n = 2i \sin(nA_n\pi)$,
we perform the first transformation of our Riemann-Hilbert problem.
We define for $z \in \mathbb C \setminus \Sigma^{(n)}$,
\begin{equation} \label{eq51}
    U(z) = (c_n)^{- \frac{1}{2} \sigma_3} e^{- \frac{1}{2} n \ell_n  \sigma_3} Y(z)
    e^{-ng_n(z) \sigma_3} e^{\frac{1}{2} n \ell_n \sigma_3} (c_n)^{\frac{1}{2} \sigma_3}.
\end{equation}
Here $\sigma_3$ is the Pauli matrix $\sigma_3 = \begin{pmatrix} 1 & 0\\0 & -1\end{pmatrix}$.

From the Riemann-Hilbert problem for $Y$
and the jump relations satisfied by $g_n$, it follows by a straightforward
calculation that $U$ is the unique solution of the following Riemann-Hilbert problem.

\subsubsection*{Riemann-Hilbert problem for $U$:}
The problem is to determine a $2\times 2$ matrix valued function
$U : \mathbb C \setminus \Sigma^{(n)} \to \mathbb C^{2\times 2}$ such that
\begin{enumerate}
\item[(a)] $U(z)$ is analytic for $z \in \mathbb C \setminus \Sigma^{(n)}$,
\item[(b)] $U(z)$ possesses continuous boundary values for $z \in \Sigma^{(n)}$,
denoted by $U_+(z)$ and $U_-(z)$, and $U_+(z) = U_-(z) V_U(z)$ where
\begin{align}
    \label{jumpVU1}
    V_U(z) &= \left( \begin{array}{cc}
    e^{2n\phi_{n,+}(z) +2nA_n\pi i} & 1  \\
    0 & e^{2n \phi_{n,-}(z) - 2nA_n\pi i} \end{array} \right)
    && \mbox{for } z \in (\beta_{1,n}, \beta_{2,n}], \\[10pt]
    \label{jumpVU2}
    V_U(z) &= \left(\begin{array}{cc}
    1 & e^{-2n\tilde{\phi}_n(z)} \\
    0 & 1 \end{array} \right)
    && \mbox{for } z \in (\beta_{2,n}, \infty), \\[10pt]
    \label{jumpVU3}
    V_U(z) &=
    \left(\begin{array}{cc}
    e^{2n \phi_n(z) \pm nA_n \pi i} & c_n^{-1} \\
    0 & e^{-2n\phi_n(z) \mp nA_n \pi i} \end{array} \right)
    && \mbox{for } z \in \Gamma_0^{(n)} \cap \mathbb C^{\pm}.
\end{align}
\item[(c)] $\ds U(z) = I + O\left(\frac{1}{z}\right)$
    as $z \to \infty$.
\end{enumerate}

\subsection{Second transformation $U \mapsto T$}

We choose a small $\varepsilon > 0$, so that the two disks
\begin{equation} \label{eqdisks}
    \Delta_{\varepsilon}(\beta_j) = \{ z \in \mathbb C \mid
    |z-\beta_j| < \varepsilon \}, \qquad j = 1,2,
\end{equation}
are disjoint.
Then there exist $n_0$ such that for every $n \geq n_0$ and $j=1,2$, we have that
$\beta_{j,n} \in \Delta_{\varepsilon}(\beta_j)$ and that the contour $\Sigma^{(n)}$
is within an $\varepsilon$-neighborhood of $\Sigma = \Gamma_0 \cup [\beta_1, \beta_2]$.

\begin{figure}[th!]
\centerline{\includegraphics[width=10cm]{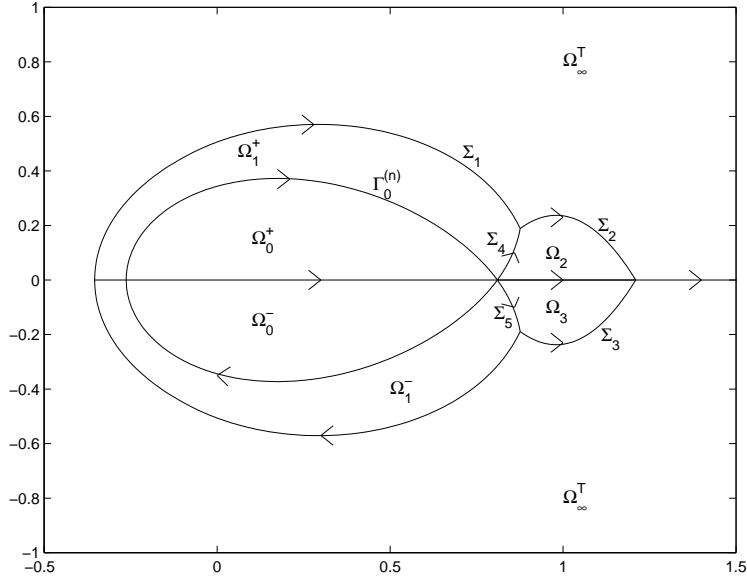}}
\caption{Contour for the Riemann-Hilbert problem for $T$,
for the value $A_n = 0.99$. The contour consists of a semi-infinite real interval
$[x_0,\infty)$,
the closed curve $\Gamma_0$, and the five arcs $\Sigma_j$, $j=1,\ldots,5$.
The contour bounds the regions
$\Omega_0^{\pm}$, $\Omega_1^{\pm}$, $\Omega_2$, $\Omega_3$, and $\Omega_{\infty}^{T}$.}
\label{fig:plotT2}
\end{figure}

\begin{figure}[th!]
\centerline{\includegraphics[width=10cm]{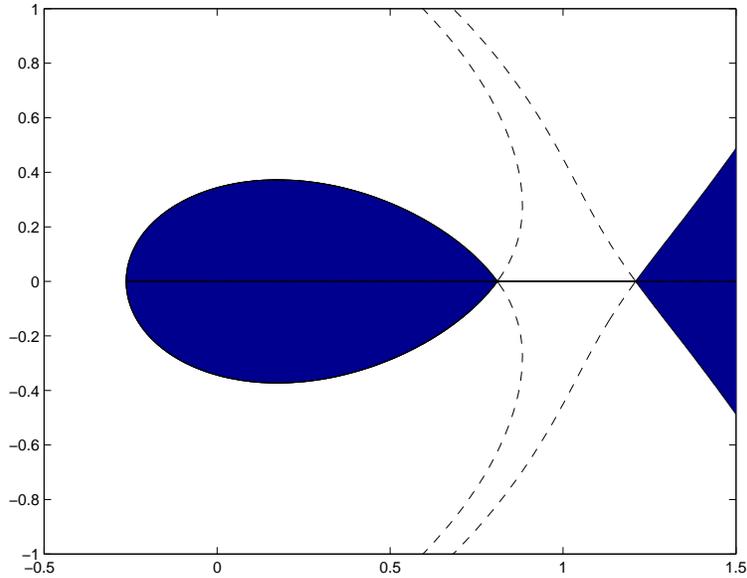}}
\caption{The dark region is the region where $\Re \phi > 0$, the white region
is the region where $\Re \phi < 0$. The parts of the contours $\Sigma_j$
for $j=1,\ldots, 5$ that are outside the disks
are chosen in the white region and within distance $2\varepsilon$ from $\Sigma_0$.
The particular plot is for the value $A = 0.99$.}
\label{fig:shaded}
\end{figure}

Now for $n \geq n_0$, we open up regions around $\Gamma_0^{(n)}$ and
$[\beta_1^{(n)},\beta_2^{(n)}]$. The regions are called $\Omega_0^{\pm}$,
$\Omega_1^{\pm}$, $\Omega_2$, and $\Omega_3$ as shown in Figure \ref{fig:plotT2}.
The exterior region is called $\Omega_{\infty}^T$.
The boundaries of the regions form the contour
$[x_0,\infty) \cup \bigcup_{j=0}^5 \Sigma_j$,
where $x_0$ is the intersection of $\Sigma_1$ with the negative real
axis. The contour is oriented as in Figure \ref{fig:plotT2}.
We can (and do) take the contours $\Sigma_j$ to be {\em independent of $n$ outside
of the disks $\Delta_{\varepsilon}(\beta_j)$} and so that $\Re \phi < 0$ on the
contours $\Sigma_j$ outside the disks.
We can clearly arrange for that. Indeed the solid lines in Figure \ref{fig:traject}
are the curves where $\Re \phi = 0$, and this includes the closed contour $\Gamma_0$.
The region where $\Re \phi < 0$ is the unbounded white region shown in
Figure \ref{fig:shaded}. The dark region is the region where
$\Re \phi > 0$. So away from the points $\beta_1$ and $\beta_2$,
we choose the contours $\Sigma_j$ in the white region.

We also take $\Sigma_j$ for $j=1,\ldots,5$ to be within distance
$\leq 2\varepsilon$ from $\Sigma_0$, and such
that the point where $\Sigma_1$, $\Sigma_2$, and $\Sigma_4$
meet, and the point where $\Sigma_1$, $\Sigma_3$, and $\Sigma_5$ meet
are outside the disks $\Delta_{\varepsilon}(\beta_j)$.

With this choice of contours we define $T$ for $n \geq n_0$, by
\begin{align} \label{eq63}
    T(z) &= U(z)  \left(\begin{array}{cc} 0 & c_n^{-1} \\
    -c_n & e^{-2n \phi_n(z) \mp nA_n\pi i} \end{array} \right)
    && \mbox{for } z \in \Omega_0^{\pm}, \\[10pt]
    \label{eq64}
    T(z) &= U(z)  \left(\begin{array}{cc} 1 & 0 \\
    -c_n e^{2n \phi_n(z) \pm nA_n\pi i} & 1 \end{array} \right)
    && \mbox{for } z \in \Omega_1^{\pm}, \\[10pt]
    \label{eq65}
    T(z) &= U(z)  \left(\begin{array}{cc} 1 & 0 \\
    -e^{2n \phi_n(z) +2 nA_n\pi i} & 1 \end{array} \right)
    && \mbox{for } z \in \Omega_2, \\[10pt]
    \label{eq66}
    T(z) &= U(z)  \left(\begin{array}{cc} 1 & 0 \\
    e^{2n \phi_n(z) -2 nA_n\pi i} & 1 \end{array} \right)
    && \mbox{for } z \in \Omega_3, \\[10pt]
    \label{eq67}
    T(z) &= U(z)
    && \mbox{for } z \in \Omega_{\infty}^T.
\end{align}

In (\ref{eq64}) we have defined $T$ with a different formula in
$\Omega_1^+$ and $\Omega_1^-$. So one might expect a jump across
$\Omega_1 \cap \mathbb R$. However, we have as in (\ref{jumpphi})
\begin{equation} \label{jumpphin}
    \phi_{n,+}(x) = \phi_{n,-}(x) - A_n \pi i, \qquad
    \mbox{for } x < 0,
\end{equation}
so that (\ref{eq64}) gives that $T_+(x) = T_-(x)$ for $x \in \Omega_1 \cap \mathbb R$.
Similarly, by (\ref{eq63}) and (\ref{jumpphin}), we have $T_+(x) = T_-(x)$
for $x \in \Omega_0 \cap \mathbb R^-$. It follows that
$T$ has an analytic continuation across $(x_0, 0)$.

\begin{figure}[th!]
\centerline{\includegraphics[width=10cm]{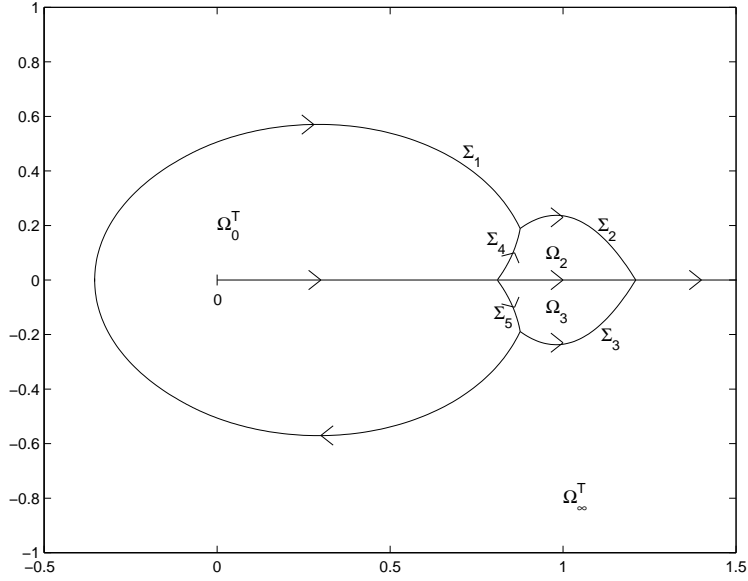}}
\caption{Reduced contour $\Sigma^T$ for the Riemann-Hilbert problem for
$T$, for the value $A_n = 0.99$. The contour consists of the semi
axis $[0,\infty)$ and the five arcs $\Sigma_j$, $j=1,\ldots,5$.
The contour bounds the regions $\Omega_0^T$, $\Omega_2$,
$\Omega_3$, and $\Omega_{\infty}^{T}$.} \label{fig:plotT3}
\end{figure}

Next, we see by direct calculation from (\ref{jumpVU3}) and (\ref{eq63})
that $T_+ = T_-$ on $\Gamma_0^{(n)}$, and therefore $T$
has an analytic continuation across $\Gamma_0^{(n)} $ as well. So $T$ is
analytic on the complement of the reduced contour $\Sigma^T$ shown
in Figure \ref{fig:plotT3}. The regions bounded by the contour
$\Sigma^T$ are denoted by $\Omega_0^T$, $\Omega_2$, $\Omega_3$,
and $\Omega_{\infty}^T$ as indicated in Figure \ref{fig:plotT3}.

Now from the Riemann-Hilbert problem satisfied by $U$, we obtain  by straightforward
calculations the jumps for $T$ on $\Sigma^T$.
The result is that $T$ is the unique solution of the following Riemann-Hilbert problem.
\subsubsection*{Riemann-Hilbert problem for $T$:}
The problem is to determine a $2\times 2$ matrix valued function
$T : \mathbb C \setminus \Sigma^T \to \mathbb C^{2\times 2}$ such that
the following hold:
\begin{enumerate}
\item[(a)] $T(z)$ is analytic for $z \in \mathbb C \setminus \Sigma^T$,
\item[(b)] $T(z)$ possesses continuous boundary values for $z \in \Sigma^T$,
denoted by $T_+(z)$ and $T_-(z)$, such that  $T_+(z) = T_-(z) V_T(z)$ where
\begin{align} \label{eq68}
    V_T(z) &= \left(\begin{array}{cc} 1 & 0 \\
    c_n e^{2n \phi_n(z) \pm nA_n\pi i} & 1 \end{array} \right)
    && \mbox{for } z \in \Sigma_1 \cap \mathbb C^{\pm}, \\[10pt]
     \label{eq69}
    V_T(z) &= \left( \begin{array}{cc}
     1 & 0 \\
     e^{2n \phi_n(z) +2nA_n\pi i} & 1 \end{array} \right)
     && \mbox{for } z\in \Sigma_2, \\[10pt]
    \label{eq610}
    V_T(z) &= \left(\begin{array}{cc} 1 & 0 \\
    e^{2n \phi_n(z) - 2nA_n\pi i} & 1 \end{array} \right)
    && \mbox{for } z \in \Sigma_3, \\[10pt]
    \label{eq611}
    V_T(z) &= \left(\begin{array}{cc} 1 & 0 \\
    e^{2n \phi_n(z)} & 1 \end{array} \right)
    && \mbox{for } z \in \Sigma_4 \cup \Sigma_5, \\[10pt]
    \label{eq612}
    V_T(z) &= \left(\begin{array}{cc} 1 & e^{-2n\phi_n(z)} \\
    0 & 1 \end{array} \right)
    && \mbox{for } z \in [0,\beta_{1,n}), \\[10pt]
    \label{eq613}
    V_T(z) &=  \left( \begin{array}{cc}
     0 & 1 \\ -1 & 0 \end{array} \right)
     && \mbox{for } z\in (\beta_{1,n}, \beta_{2,n}), \\[10pt]
    \label{eq614}
    V_T(z) &= \left( \begin{array}{cc}
     1 & e^{-2n \tilde{\phi}_n(z)} \\
     0 & 1 \end{array} \right)
     && \mbox{for } z\in (\beta_{2,n}, \infty).
\end{align}
\item[(c)] $\ds    T(z) = I + O\left(\frac{1}{z}\right)$
     as $ z \to \infty$.
\end{enumerate}

\begin{remark}
The contour $\Sigma^T$ has been constructed so that
$\Sigma_j \setminus (\Delta_{\varepsilon}(\beta_1) \cup \Delta_{\varepsilon}(\beta_2))$
is $n$-independent, and so that $\phi(z) < 0$ there. Since $\phi_n(z) \to \phi(z)$
as $n \to \infty$, we see that the jump matrices in (\ref{eq68})--(\ref{eq611})
converge to the identity matrix as $n \to \infty$. The convergence is uniform for
$z$ outside the disks $\Delta_{\varepsilon}(\beta_j)$.

Similarly, since $\phi(z) > 0$ for $z \in [0,\beta_{1})$,
$\tilde{\phi}(z) > 0$ for $z \in (\beta_{2}, \infty)$, and $\phi_n \to \phi$,
and $\tilde{\phi}_n \to \tilde{\phi}$ as $n \to \infty$, the
jump matrices in (\ref{eq612}) and (\ref{eq614}) also converge to the identity
matrix as $n \to \infty$, and the convergence is uniform for $z$ outside the
disks $\Delta_{\varepsilon}(\beta_j)$.
\end{remark}

\subsection{Construction of the parametrix for $T$}

Observe that we may write the jump matrices in (\ref{eq69}) and (\ref{eq610}) also as
$\left(\begin{array}{cc} 1 & 0 \\ e^{2n\tilde{\phi}_n(z)} & 1 \end{array} \right)$.
Then it follows that we have a situation which is basically the same as the one
considered by Deift in \cite[p.~199]{Deift}. There the part $\Sigma_1$ is absent,
and the contour extends to $-\infty$ along the negative real axis. However, these
differences are not important since on these parts
the jump matrices are uniformly $I + O(e^{-cn|z|})$ as $n \to \infty$, for some $c > 0$.

It follows that we can do the same analysis as in \cite{Deift} (see also \cite{DKMVZ1,DKMVZ2})
to construct the parametrix for $T$. The behavior away from $\beta_1$ and $\beta_2$
is governed by the solution of the following Riemann-Hilbert problem.

\subsubsection*{Riemann-Hilbert problem for $N$:}
The problem is to determine $N : \mathbb C \setminus [\beta_{1,n},\beta_{2,n}]
\to \mathbb C^{2\times 2}$
such that the following hold.
\begin{enumerate}
\item[(a)] $N(z)$ is analytic for $z \in \mathbb C \setminus [\beta_{1,n},\beta_{2,n}]$,
\item[(b)] $N(z)$ possesses continuous boundary values for $z \in
(\beta_{1,n}, \beta_{2,n})$, denoted by $N_+(z)$ and $N_-(z)$, such that
\begin{equation} \label{eq71}
    N_+(z) = N_-(z)
    \left(\begin{array}{cc} 0 & 1 \\ -1 & 0 \end{array} \right)
    \qquad \mbox{ for } z \in (\beta_{1,n}, \beta_{2,n}),
\end{equation}
and $N$ has at most $\frac{1}{4}$-root singularities at $\beta_{1,n}$ and $\beta_{2,n}$,
\item[(c)] $N(z) = I + O\left(\frac{1}{z}\right)$ as $z \to \infty$.
\end{enumerate}

Writing $a_n(z) = \frac{(z-\beta_{2,n})^{1/4}}{(z-\beta_{1,n})^{1/4}}$,
we then have that the solution of the Riemann-Hilbert problem for $N$ is
\begin{equation} \label{eq72}
    N(z) =
    \left(\begin{array}{cc}
    \frac{a_n(z) + a_n(z)^{-1}}{2} &
    \frac{a_n(z) - a_n(z)^{-1}}{2i} \\[10pt]
    \frac{a_n(z) - a_n(z)^{-1}}{-2i} &
    \frac{a_n(z) + a_n(z)^{-1}}{2}
    \end{array} \right)
    \qquad \mbox{ for } z \in \mathbb C \setminus [\beta_{1,n}, \beta_{2,n}].
\end{equation}
\medskip

Near $\beta_1$ and $\beta_2$ we construct local parametrices, denoted by
$P_1$ and $P_2$, respectively.
Recall that we have chosen $\varepsilon$ so that the disks $\Delta_{\varepsilon}(\beta_1)$
and $\Delta_{\varepsilon}(\beta_2)$ are disjoint.

\subsubsection*{Riemann-Hilbert problem for $P_{j}$, $j=1,2$:}
The problem is to determine matrix valued functions
$P_j : \Delta_{\varepsilon}(\beta_j) \setminus \Sigma^T \to \mathbb C^{2\times 2}$
such that for $j=1,2$,
\begin{enumerate}
\item[(a)] $P_j(z)$ is analytic for $z \in \Delta_{\varepsilon}(\beta_j) \setminus \Sigma^T$,
\item[(b)] $P_j$ satisfies the same jump conditions on $\Sigma^T \cap \Delta_{\varepsilon}(\beta_j)$
as $T$ does.
\item[(c)]
We have as $n \to \infty$,
\begin{equation} \label{eq73}
    P_j(z) = \left(I + O \left(\frac{1}{n} \right) \right) N(z)
    \qquad \mbox{uniformly for } |z - \beta_j| = \varepsilon.
\end{equation}
\end{enumerate}

The Riemann-Hilbert problems for $P_1$ and $P_2$ are solved as
in \cite{Deift,DKMVZ2} with the use of Airy functions, see also
\cite{DKMVZ1,KuMc}. Note that $P_1$ and $P_2$,
as well as $N$, depend on $n$.

The formula for $P_2$ is as in \cite[Eq.\ (7.24)]{DKMVZ2}.
It uses a conformal mapping
\[ \zeta = f_n(z) = \left( \frac{3}{2} n \tilde{\phi}_n(z) \right)^{3/2} \]
from $\Delta_{\varepsilon}(\beta_2)$ onto a neighborhood of $0$ in the $\zeta$-plane,
and an explicit $2\times 2$ matrix $\Psi^{\sigma}(\zeta)$ built
out of Airy functions which is given in \cite[p.~1522]{DKMVZ2}. The superscript
$\sigma$ denotes an angle that may be any number in $(\pi/3, \pi)$, possibly
be $n$-dependent. Then $P_2$ has the form
\[ P_2(z) = E_n(z) \Psi^{\sigma}(f_n(z))
    \begin{pmatrix} e^{n\tilde{\phi}_n(z)} & 0 \\
    0 & e^{-n\tilde{\phi}_n(z)} \end{pmatrix} \]
where $E_n$ is an analytic pre-factor that takes care of the matching
condition (\ref{eq73}), see \cite[Eq.\ (7.23)]{DKMVZ2} for the precise form
of $E_n$.
The construction of the parametrix $P_1$ near the left end-point $\beta_1$
is similar, see \cite[pp.\ 1526-1527]{DKMVZ2}.

\subsection{Final transformation $T \mapsto S$}

Having $N$, $P_1$, and $P_2$, we define
\begin{eqnarray} \label{eq81}
    S(z) & = & T(z) N(z)^{-1} \qquad
    \mbox{for } z \in \mathbb C \setminus
        (\Sigma^T \cup \overline{\Delta_{\varepsilon}(\beta_1)}
            \cup \overline{\Delta_{\varepsilon}(\beta_2)}), \\[10pt]
    \label{eq82}
    S(z) & = & T(z) P_1(z)^{-1} \qquad
    \mbox{for } z \in \Delta_{\varepsilon}(\beta_1) \setminus \Sigma^T, \\[10pt]
    \label{eq83}
    S(z) & = & T(z) P_2(z)^{-1} \qquad
    \mbox{for } z \in \Delta_{\varepsilon}(\beta_2) \setminus \Sigma^T.
\end{eqnarray}
\begin{figure}[th!]
\centerline{\includegraphics*[width=10cm,height=8cm]{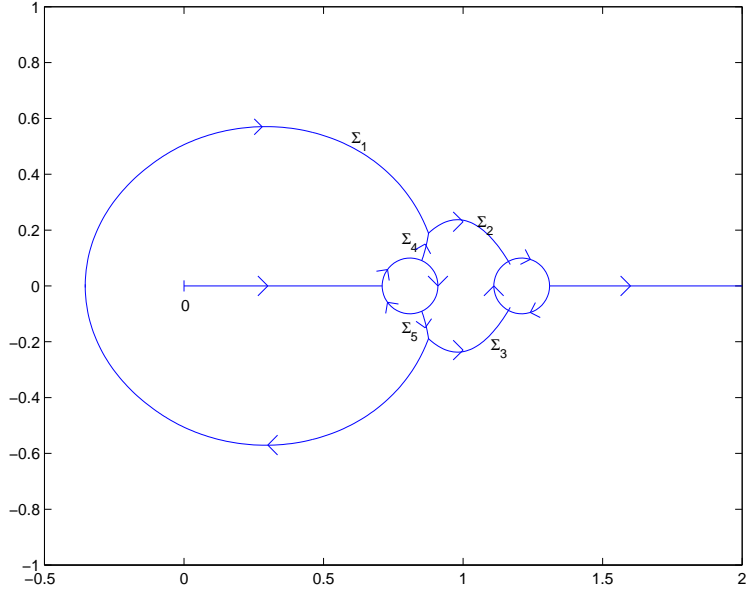}}
\caption{Contour $\Sigma^S$ for the Riemann-Hilbert problem for $S$,
for the value $A = 0.99$. The two circles are the boundaries of two
$\varepsilon$-neighborhoods of $\beta_1$ and $\beta_2$, with $\varepsilon = 0.1$.}
\label{fig:plotS}
\end{figure}
Then $S$ is defined and analytic on $\mathbb C \setminus (\Sigma^T
\cup \partial \Delta_{\varepsilon}(\beta_1) \cup \partial
\Delta_{\varepsilon}(\beta_2))$. However it follows from the
construction that $S$ has the identity jump on the interval
$(\beta_1+\varepsilon,\beta_2-\varepsilon)$ and on the parts of $\Sigma^T$
that are within the disks $\Delta_{\varepsilon}(\beta_j)$.
Therefore $S$ has an analytic continuation to $\mathbb C \setminus
\Sigma^S$, where $\Sigma^S$ is the system of contours shown in
Figure \ref{fig:plotS}. The contour $\Sigma^S$ does not depend on $n$.

Then $S$ satisfies the following Riemann-Hilbert problem.
\subsubsection*{Riemann-Hilbert problem for $S$:}
The problem is to determine $S : \mathbb C \setminus \Sigma^S \to \mathbb C^{2\times 2}$
such that the following hold.
\begin{enumerate}
\item[(a)] $S(z)$ is analytic for $z \in \mathbb C \setminus \Sigma^S$,
\item[(b)] $S(z)$ possesses continuous boundary values for $z \in \Sigma^S$,
denoted by $S_+(z)$ and $S_-(z)$, such that $S_+(z) = S_-(z) V_S(z)$ where
\begin{align} \label{eq84}
    V_S(z) &= P_1(z) N(z)^{-1}
    && \mbox{for } \partial U_{\varepsilon}(\beta_1), \\[10pt]
     \label{eq85}
    V_S(z) &= P_2(z) N(z)^{-1}
    && \mbox{for } \partial U_{\varepsilon}(\beta_2), \\[10pt]
     \label{eq86}
    V_S(z) &= N(z) V_T(z) N(z)^{-1}
    && \mbox{for } z \in \Sigma^S \setminus
        \left(\partial U_{\varepsilon}(\beta_1)
        \cup \partial \tilde{U}_{\varepsilon}(\beta_2)\right)
\end{align}
\item[(c)]  $\ds S(z) = I + O\left(\frac{1}{z}\right)$ as $z \to \infty$.
\end{enumerate}

All jump matrices for $S$ are close to the identity matrix if $n$ is large. The jump matrices
in (\ref{eq84}) and (\ref{eq85}) are $I + O(1/n)$, and the jump matrices in (\ref{eq86})
are uniformly $I + O(e^{-cn})$ for some $c > 0$.
In addition, the jump matrix (\ref{eq86}) on the unbounded interval
$(\beta_2 + \varepsilon, \infty)$ is uniformly $I + O(e^{-cn|z|})$
for some $c > 0$, so that the jump matrix is close to the identity matrix
both in $L_2$-norm and in $L_{\infty}$-norm on the contour $\Sigma^S$.
Also note that by construction the contour $\Sigma^S$ does not depend on $n$.
Then it follows as in \cite{Deift} that
\begin{equation} \label{eq811}
S(z) = I + O\left(\frac{1}{n}\right)
\qquad \mbox{ uniformly for } z \in \mathbb C \setminus \Sigma^S,
\end{equation}
as $n \to \infty$.

\subsection{Proof of Theorem 1.2}
Tracing back the steps $Y \mapsto U \mapsto T \mapsto S$ and
using (\ref{eq811}) we obtain strong asymptotics for $Y$ in every
region in the complex plane. In particular since
\[ Y_{11}(z) = P_n(z) = \frac{n!}{(-n)^n} L_n^{(-n A_n)}(nz), \]
this yields strong asymptotics for the generalized Laguerre
polynomials $L_n^{(-n A_n)}(nz)$.

The formulas we obtain in various regions are similar to the
ones obtained in \cite{DKMVZ2} for orthogonal polynomials with
respect to exponential weights on $\mathbb R$.
We will not give all these results here, but restrict ourselves
to the ones that are needed for the proof of Theorem 1.2.

\begin{varproof} {\bf of Theorem 1.2.}
As noted before, we may and do assume that $\alpha_n \not\in \mathbb Z$
for every $n$. For every $n$, we can then set up the Riemann-Hilbert problem
for $Y$, and by  (\ref{solutionY}) and (\ref{eq51}) we have
\begin{equation} \label{eq91}
    P_n(z) = Y_{11}(z) = U_{11}(z) e^{ng_n(z)}.
\end{equation}
Using the formulas (\ref{eq63})--(\ref{eq67}) which express
$T$ in terms of $U$, we find $P_n$ in terms of the entries of $T$.
The precise formula depends on the region.

For $z$ in the outer region $\Omega_{\infty}^T$, we find by (\ref{eq67}) and (\ref{eq91}),
\begin{equation} \label{eq92}
    P_n(z) = T_{11}(z) e^{ng_n(z)}.
\end{equation}
If, in addition, $|z-\beta_2| > \varepsilon$, then $T = SN$ by (\ref{eq81})
and $S = I + O(1/n)$ by (\ref{eq811}), so that
\begin{equation} \label{eq93}
    P_n(z) e^{-n g_n(z)} = N^{(n)}_{11}(z) \left(1 + O\left(\frac{1}{n}\right)\right),
\end{equation}
uniformly for $z \in \Omega_{\infty}^T \setminus \Delta_{\varepsilon}(\beta_2)$.
We have written $N^{(n)}$ to emphasize the $n$-dependence of $N$.
As $n\to \infty$, we have by (\ref{eq72})
\begin{equation} \label{aisliman}
    \lim_{n\to\infty} N^{(n)}_{11}(z) =
    \lim_{n\to\infty} \frac{a_n(z) + a_n(z)^{-1}}{2} =
    \frac{a(z) + a(z)^{-1}}{2}
\end{equation}
uniformly for $z \in \Omega_{\infty}^T \setminus \Delta_{\varepsilon}(\beta_2)$,
where $a(z) = \frac{(z-\beta_2)^{1/4}}{(z-\beta_1)^{1/4}}$. Since
the right-hand side of (\ref{aisliman}) does not have zeros in $\mathbb C$,
it follows from (\ref{eq93}), (\ref{aisliman}), and Hurwitz' theorem,
that there are no zeros of $P_n$ in
$\Omega_{\infty}^T \setminus \Delta_{\varepsilon}(\beta_2)$,
if $n$ is sufficiently large. From a contour deformation argument as in
\cite[p.1532]{DKMVZ2}, we get that (\ref{eq93}) holds for
$z \in (\beta_2 + \varepsilon, \infty)$ as well. Thus, by the same argument,
there are no zeros in this interval if $n$ is sufficiently large, and so,
for large $n$, all zeros of $P_n$ are in
$\overline{\Omega_0^T \cup \Omega_2 \cup \Omega_3 \cup \Delta_{\varepsilon}(\beta_2)}$.
\medskip

Next, we have by (\ref{eq63}) and (\ref{eq91})
\begin{equation} \label{PninOmega0}
    P_n(z) = U_{11}(z)e^{ng_n(z)} =
    \left[ e^{-2n\phi_n(z) \mp nA_n\pi i} T_{11}(z) + c_n T_{12}(z)
    \right] e^{ng_n(z)}
    \qquad \mbox{for } z \in \Omega_0^{\pm},
\end{equation}
and by (\ref{eq64}) and (\ref{eq91})
\begin{equation} \label{PninOmega1}
    P_n(z) = U_{11}(z)e^{ng_n(z)} =
    \left[T_{11}(z) + c_n e^{2n\phi_n(z) \pm nA_n\pi i} T_{12}(z)
    \right] e^{ng_n(z)}
    \qquad \mbox{for } z \in \Omega_1^{\pm}.
\end{equation}
For $|z-\beta_j| > \varepsilon$, we have $T = SN$ by (\ref{eq81})
and $S = I+ O(1/n)$ by (\ref{eq811}). Then we get from
(\ref{PninOmega0})
\begin{eqnarray} \nonumber
    P_n(z) e^{-ng_n(z)} & = &
    \left[ e^{-2n\phi_n(z) \mp nA_n\pi i} N_{11}^{(n)}(z)
        +c_n N_{12}^{(n)}(z) \right]
        \left(1 + O\left(\frac{1}{n}\right) \right) \\[10pt]
    & & + \label{asympPninOmega0}
    \left[ e^{-2n\phi_n(z) \mp nA_n\pi i} N_{21}^{(n)}(z)
         +c_n N_{22}^{(n)}(z) \right]
        O\left(\frac{1}{n}\right),
\end{eqnarray}
uniformly for $z \in \overline{\Omega_0} \setminus \Delta_{\varepsilon}(\beta_1)$,
and from (\ref{PninOmega1})
\begin{eqnarray} \nonumber
    P_n(z) e^{-ng_n(z)} & = &
    \left[ N_{11}^{(n)}(z) + c_n e^{2n\phi_n(z) \pm nA_n\pi i}
        N_{12}^{(n)}(z) \right]
        \left(1 + O\left(\frac{1}{n}\right) \right) \\[10pt]
    & & + \label{asympPninOmega1}
    \left[ N_{21}^{(n)}(z) + c_n e^{2n\phi_n(z) \pm nA_n\pi i}
        N_{22}^{(n)}(z) \right]
        O\left(\frac{1}{n}\right),
\end{eqnarray}
uniformly for $z \in \overline{\Omega_1} \setminus \Delta_{\varepsilon}(\beta_1)$.

As $n \to \infty$, we clearly have
\[ \lim_{n \to \infty}  N^{(n)}(z) =
    \left(\begin{array}{cc}
    \frac{a(z) + a(z)^{-1}}{2} &
    \frac{a(z) - a(z)^{-1}}{2i} \\[10pt]
    \frac{a(z) - a(z)^{-1}}{-2i} &
    \frac{a(z) + a(z)^{-1}}{2}
    \end{array} \right), \]
from which it follows that the entries of $N^{(n)}$ are uniformly
bounded and uniformly bounded away from zero on any compact set which
does not contain $\beta_1$ and $\beta_2$.
We also recall that
\begin{equation} \label{Lislimcn}
    L = \lim_{n \to \infty} |c_n|^{1/n}.
\end{equation}
Now we distinguish the cases $L=0$ and $L= e^{-r}$ as in Theorem 1.2.

\begin{itemize}
\item[(a)] Case $L = 0$.

Let $K = \overline{\Omega_0^T} \setminus
    \left( \Delta_{\varepsilon}(\beta_1) \cup \Delta_{\varepsilon}(0) \right)$.
Then the functions $\phi_n$ are uniformly bounded on $K$,
so that we get from (\ref{Lislimcn}) and $L = 0$ that
\[ \lim_{n \to \infty} \left| c_n e^{2n\phi_n(z) \mp A_n \pi i} \right|^{1/n} = 0
    \qquad \mbox{uniformly for } z \in K.
\]
Since the entries of $N^{(n)}$ are uniformly bounded and  uniformly bounded
away from $0$, it then follows from (\ref{asympPninOmega0})
and (\ref{asympPninOmega1}) that
\[ P_n(z) e^{-ng_n(z)} =
    e^{-2n\phi_n(z) \mp nA_n\pi i} \frac{a_n(z) + a_n(z)^{-1}}{2}
        \left(1 + O\left(\frac{1}{n}\right) \right) \]
uniformly for $z \in K \cap \overline{\Omega_0}$, and
\[ P_n(z) e^{-ng_n(z)} =
 \frac{a_n(z) + a_n(z)^{-1}}{2} \left(1 + O\left(\frac{1}{n}\right) \right)
 \]
uniformly for $z \in K \cap \overline{\Omega_1}$.
Thus for large $n$, there are no zeros in $K$, so that
all zeros are in $\Delta_{\varepsilon}(0)
\cup \Delta_{\varepsilon}(\beta_1) \cup \overline{\Omega_2}
\cup \overline{\Omega_3} \cup \Delta_{\varepsilon}(\beta_2)$.
Since $\varepsilon$ can be chosen arbitrarily close to $0$, and the
domains $\Omega_2$ and $\Omega_3$ lie within the $2\varepsilon$ neighborhood
of $[\beta_1,\beta_2]$,  it follows that
the zeros accumulate only on $\{0\} \cup [\beta_1,\beta_2]$.

\item[(b)] Case $L=e^{-r}$ with $0 \leq r < \infty$.

To locate zeros of $P_n$ we have to balance two contributions in
(\ref{asympPninOmega0}) or (\ref{asympPninOmega1}),
which can only be done near the curve $\Gamma_r$ where $\Re \phi = r/2$. Indeed,
if $K$ is
\[ K = \{ z \in \overline{\Omega_0^T} \mid \dist(z, \Gamma_r) \geq \varepsilon \} \]
then there is $\delta > 0$ such that
$|\phi_n(z) - r/2| > \delta$ for $z \in K$ and for $n$ large enough.
Then either $c_n N_{12}^{(n)}(z) $ or $e^{-2n\phi_n(z) \mp n A_n\pi i} N_{11}^{(n)}(z)$
dominates the other term in the asymptotic formula (\ref{asympPninOmega0}).
The former dominates if $\Re \phi_n(z) > r/2 + \delta$ (which is for $z$ in the
interior region bounded by $\Gamma_r$), and the latter dominates
if $\Re \phi_n(z) < r/2 - \delta$. Then it follows that for large $n$, there are
no zeros in $K \cap \overline{\Omega_0}$. Similarly, it follows from (\ref{asympPninOmega0})
that no zeros in $K \cap \overline{\Omega_1}$ for large $n$.
Thus the zeros in $\Omega_0^T \setminus \Delta_{\varepsilon}(\beta_1)$
are in a $\varepsilon$ neighborhood of $\Gamma_r$ if $n$ is large.
In addition, there can be zeros in $\Delta_{\varepsilon}(\beta_1)
\cup \Delta_{\varepsilon}(\beta_2)
\cup \overline{\Omega_2} \cup \overline{\Omega_3}$.
Letting $\varepsilon \to 0$, we see that zeros can accumulate on $\Gamma_r$
and on the interval $[\beta_1,\beta_2]$ only.
\end{itemize}

So now we know that the zeros accumulate on $\{0\} \cup [\beta_1,\beta_2]$
in case $L=0$, and on $\Gamma_r \cup [\beta_1, \beta_2]$ in case $L = e^{-r} > 0$.
It remains to determine the asymptotic zero distribution. To that end, we note
that we get from (\ref{eq93}),
\begin{equation} \label{zerodist1}
    \lim_{n\to\infty} \left(\frac{1}{n} \log |P_n(z)| - \Re g_n(z) \right) = 0,
\end{equation}
uniformly for $z$ in compact subsets of $\Omega_{\infty}^T$.
From the definition of $g_n$, see (\ref{defgnz}), it follows that
\begin{equation} \label{zerodist2}
    \lim_{n \to \infty} \Re g_n(z) = \int \log|z-s| \, d\mu_0(s)
\end{equation}
uniformly for $z$ in compact subsets of $\mathbb C \setminus
\left(\Gamma_0 \cup [\beta_1,\beta_2]\right)$. Here $\mu_0$ is the
measure on $\Gamma_0 \cup [\beta_1, \beta_2]$ given by
(\ref{defmur}) with $r = 0$.
Thus by (\ref{zerodist1}) and (\ref{zerodist2}),
\begin{equation} \label{zerodist3}
    \lim_{n \to \infty} \frac{1}{n} \log |P_n(z)| =
   \int \log|z-s| \, d\mu_0(s),
\end{equation}
for $z$ in $\Omega_{\infty}^T$.

From the definition of the measures $\mu_r$ in (\ref{defmur}) it follows
that $\int \log|z-s| d\mu_r(s)$ does not depend on $r$ if $z \in \Omega_{\infty}$.
Indeed, the parts of $\mu_r$ on the interval $[\beta_1, \beta_2]$ are
the same for every $r$, and for the part $\nu_r$ of $\mu_r$ on $\Gamma_r$,
we have by Cauchy's theorem
\begin{equation} \label{indepr}
    \int \log(z-s) d\nu_r =
    \frac{1}{2\pi i} \int_{\Gamma_r} \log(z-s) \frac{R(s)}{s} ds
    = (\log z) R(0) = (\log z) A,
\end{equation}
for $z \in \Omega_{\infty}^T$, which is also $r$-independent.
The relation (\ref{indepr}) is valid for
any branch of  $\log(z-s)$ which is single-valued for $s$ in the
domain enclosed by $\Gamma_0$.
So $\int \log|z-s| d\mu_r(s)$ is indeed independent of $r$, so
that by (\ref{zerodist3}),
\begin{equation} \label{logPn}
 \lim_{n\to\infty} \frac{1}{n} \log |P_n(z)| =
   \int \log|z-s| \, d\mu_r(s)
   \qquad \mbox{for } z \in \Omega_{\infty}^T.
\end{equation}
Since, in case $L = e^{-r}$, all zeros accumulate on
$\Gamma_r \cup [\beta_1,\beta_2]$ and $\mu_r$ is a probability measure
on $\Gamma_r \cup [\beta_1, \beta_2]$, it follows from (\ref{logPn})
that $\mu_r$ is the asymptotic distribution of the zeros of $P_n$,
see \cite[Theorem 2.3]{MhaskarSaff} and also
\cite[Chapter III]{SaffTotik}. Note that it is important here
that $\Gamma_r \cup [\beta_1, \beta_2]$ is the boundary of its
polynomial convex hull.

In case $L=0$, it follows in a similar way that (\ref{MarPas})
is the asymptotic zero distribution.
\end{varproof}

\end{document}